\documentclass[12pt,  a4 paper]{article}
\usepackage{amsmath}
\usepackage{amsthm}
\usepackage{amssymb}
\usepackage[hidelinks]{hyperref}
\newtheorem{theorem}{Theorem}

\usepackage[left=2cm, right=2cm ]{geometry}
\date{}
\begin{document}
\begin{centering}
\title{The geometric series\\ formula and its applications}
\end{centering}
\maketitle
\author{\centerline{\large Cletus Bijalam Mbalida}
\begin{abstract}
Let $n$ be an integer and $W_n$ be the Lambert $W$ function. Let $\log$ denote the natural logarithm so that $\delta=-W_n(-\log2)/\log2$. Given that $a$ and $r$ are respectively the first term and the constant ratio of an infinite geometric series, it is proved that the limit of convergence of the geometric series is $\displaystyle\lim_{n\to\pm\infty}{a\big[r^\delta-1\big]\big[r-1\big]^{-1}}$ where $r\neq1$.\par
By applying the geometric series formula above, it is further proved that the harmonic series $\zeta(1)$ is given by $\zeta(1)=-2\big[\log2+W_n(-\log2)\big]$ and as $n\rightarrow\pm\infty$, the value of $\zeta(1)$ grows very slowly toward $\tilde\infty$, confirming the divergence of the harmonic series.
\end{abstract}
\section{Introduction}
Let $a$ and $r$ denote the first term and the constant (common) ratio of an infinite geometric series, respectively. By definition, the geometric series is generated as
\begin{align*}
a+ar+ar^2+ar^3+ar^4+\cdots.\tag{1.1}
\end{align*}
We call the sum of the above geometric series $\Pi$, without really worrying much about whether the series has a sum or not. It is already well known in the literature that if we denote the n-th term of the series by $U_n$, then as long $\big|\frac{U_{n+1}}{U_n}\big|=|r|<1$ the geometric series would converge to a finite sum of $\Pi=\frac{a}{1-r}$ and for every other value of $r$, $\Pi=\infty$.\par
In the $|r|>1$ cases, there are really no known approaches by which one can rigorously deal with the infinite geometric series of such nature. In this article, we present a novel approach by which one can deal with any infinite geometric series whose $r\neq1$.\par
In the process we are also able to deal some non-geometric infinite series and present new relations such as the expansion for $\log(x-1)$ for any real or complex number $x$ where $x\neq1$.\par
We use the following notation. The expression $\log(x)$ will always denote the natural logarithm and $n$ will be an integer. Instead of using the usual $W$ for the Lambert $W$ function, we use $W_n$ as a generalized form at each $n$. We use $p$ and $p_n$ for a prime number and the n-th prime number, respectively. For any complex variable, we use the letter $s$ and denote its real part by  $\operatorname{Re}(s)$. The letters $a$, $b$, $c$, and $r$ always represent real or complex numbers.
\begin{theorem}
Suppose $\chi$ is the finite sum of the infinite geometric series 
\begin{align*}
a+ar+ar^2+ar^3+\cdots
\end{align*}
for all $r$ except $r=1$. Then
\begin{align*}
\chi=a\frac{r^\delta-1}{r-1}, \,\, \,\, \delta=-\frac{W_n(-\log2)}{\log2}\tag{1.2}
\end{align*}
\end{theorem}
From the formula in (1.2) one sees that $\chi$ depends on $\delta$ and $\delta$ also depends on the $W_n$ . The $W_n$ function is a multivalued function whose principal value, according to the work of \cite{Corless et al}, is usually the one at $n=0$. Since every $n$ gives $\delta$ a unique solution due to the nature of $W_n$, we will denote the n-th value of $\delta$ by $\delta_n$, which will give a corresponding $\chi_n$ as the finite sum of the series. On that note we can write that the principal value of $\delta_n$ is $\delta_0=-\frac{W_0(-\log2)}{\log2}$, since the principal value of $W_n$ is usually taken at $n=0$. It is equally known that in most cases if all the values of $W_n$ are complex, then the values $W_n$ and $W_{-(n+1)}$ are complex conjugates and we will see that later in section 3.\par
It turns out when we substitute $\delta_0=-\frac{W_0(-\log2)}{\log2}$ together with known values of the parameters $a$, $r$ into the geometric series formula in (1.2), the quantity $\chi_0$--- the finite sum of the series at $n=0$ --- alone is insufficient to tell us about whether the series converges or diverges. Even if the geometric series converges, one notices that the value $\chi_0$ in most cases is not close at all to the actual sum of the series. However, as we substitute the other $\delta_n$ values into the geometric series formula together with the parameters $a$ and $r$, we observe something interesting: as $n\rightarrow\pm\infty$, $\chi_n\rightarrow\frac{a}{1-r}$ in cases where $|r|<1$, but in cases where $|r|>1$ $\chi_n\rightarrow\tilde\infty$ as $n\rightarrow\pm\infty$. The latter observation indicates that the infinite geometric series diverges for $|r|>1$. Therefore, in the $|r|>1$ cases the geometric series has infinitely many complex solutions which occur in conjugate pairs. For instance, the infinite geometric series 
\begin{align*}
\Gamma=1+\frac{1}{2}+\frac{1}{2^2}+\cdots\tag{1.3}
\end{align*}
is known from the literature to have a finite sum of $\Gamma=2$. The following are a few of the solutions of (1.3) from our new approach to dealing with infinite geometric series. Note how the solutions approach 2 as $n\rightarrow\pm\infty$. This shows the series in (1.3) indeed converges to 2.
\begin{multline*}
\Gamma_0=1.4742143439565294465800856679560310297654568506549764091... \\-
0.99933887136663452730337627926333490785340383402837588188... i\\
\Gamma_{-1}=1.4742143439565294465800856679560310297654568506549764091...\\ +
0.99933887136663452730337627926333490785340383402837588188... i\\
\Gamma_{1000}=1.9999996801764241784705239085896264398069291735087686...\\ -
0.00022058004220891148156095076445122825171717796784197993... i
\end{multline*}
\begin{multline*}
\Gamma_{-1001}=1.9999996801764241784705239085896264398069291735087686...\\ +
0.00022058004220891148156095076445122825171717796784197993... i\\
 \Gamma_{100000}=1.999999999951831540897213403424847427024870806137407...\\ -
2.206350484675363065669178309195567447968599569711646...\times 10^{-6} 
\end{multline*}
In section 3, we will learn more about  the above series as well other convergent and divergent infinite geometric series and some other non-geometric infinite series whose sums can be evaluated using the geometric series formula in (1.2) in a special way.\par
One of the non-geometric infinite series we consider in this paper is the harmonic series, usually regarded as the value of $\zeta(s)$ at $s=1$, where $\zeta(s)$ is the usual Riemann zeta function. Let $\eta(s)$ be the Dirichlet eta function. It is well known that the following relation  exists between $\zeta(s)$ and $\eta(s)$ for all $\operatorname{Re}(s)>0$.
\begin{align*}
\eta(s)=\big(1-2^{1-s}\big)\zeta(s)\tag{1.4}
\end{align*}
By applying the argument of the expansion for $\log(1+x)$ at $x=1$, one can show that $\eta(1)=\log2$. We ask whether or not it is also possible to obtain this same result in the context of (1.4) since this formula does not hold at $s=1$. In our quest for such a possiblity, it turns out (1.4) is approximated on the basis of convergence of infinite geometric series if one considers only the series $\eta(s)$ to derive (1.4) without adding or subtracting any external infinite series  to/from $\eta(s)$. By applying the formula in (1.2) it is proved that $\zeta(s)$ and $\eta(s)$ are exactly connected by
\begin{align*}
\eta(s)=\bigg[\frac{2^{s(\delta+1)}-2(2^{\delta s})+1}{2^{\delta s+s}-2^s}\bigg]\zeta(s)\tag{1.5}
\end{align*}
This relation then allows us to show that at $s=1$, $\zeta(1)=-2\big[\log2+W_n(-\log2)\big]$, and as $n\rightarrow\pm\infty$ the term $\zeta(1)\rightarrow\tilde\infty$, indicating that the harmonic series indeed diverges. We give two other new representations for the harmonic series in \S \S 4.3.1.\par
In most cases, the irrational constant $\pi$ shows up in the context of certain infnite series. In \cite{marco}, García and Marco prove that the super-regularized product over all $p$ is $4\pi^2$. Applying their result, we prove that $\pi$ and the product over all $p_n-1$ are related by
\begin{align*}
2\pi^2=\bigg(\log2\sum_{n=0}^{\infty}{2^n}\bigg)\bigg(\prod_{p_n}{\big(p_n-1\big)}\bigg)\tag{1.6}
\end{align*}
\par Furthermore, in terms of an infinite exponentiation we prove that the infinite geometric series $\sum_{n=0}^{\infty}{2^n}$ can be expressed as 
\begin{align*}
1+\sum_{n=0}^{\infty}{2^n}=2^{\scriptscriptstyle 2^{2^{2^{2^{\cdots}}}}}\tag{1.7}
\end{align*}
from which one sees that $\sum_{n=0}^{\infty}{2^n}=-1-W_n(-\log2)/\log2$. We then further show that the transcendental number $\log2$ can also be expressed as
\begin{align*}
\log2=\frac{\log(1+\sum_{n=0}^{\infty}{2^n})}{1+\sum_{n=0}^{\infty}{2^n}}\tag{1.8}
\end{align*}
\par For $\operatorname{Re}(s)>1$, each factor of the Euler product formula for $\zeta(s)$ is $\frac{p^s}{p^s-1}$ and that is equivalent to the infinite geometric series $\sum_{n=0}^{\infty}{p^{-ns}}$. Since $\frac{p^s}{p^s-1}$ is written using the geometric series formula $\frac{a}{1-r}$, there is a small error term associated with $\frac{p^s}{p^s-1}$. What is the magnitude of this error term that is associated with each $\frac{p^s}{p^s-1}$ factor of the Euler product formula for $\zeta(s)$? It is known that these error terms are insignificant but how do we prove that this is indeed the case? Let $\Delta(p; s)$ be the magnitude of the error term associated with each $\frac{p^s}{p^s-1}$ factor at a certain $s$. It is proved that $\Delta(p;s)$ is given by
\begin{align*}
\Delta(p;s)=-\frac{p^{s+sW_n(-\log2)/\log2}}{p^s-1}
\end{align*}
from which one sees that since $\operatorname{Re}(s)>1$, then as $n\rightarrow\pm\infty$ the term $\Delta(p;s)\rightarrow0$ very rapidly, indicating that indeed the error term associated with each factor of the Euler product formula for $\zeta(s)$ does not really mean anything.\par
\section{The proof of Theorem 1}
\begin{proof}
An infinite geometric series with a finite sum $\Pi$, a first term $a$ and a constant $r$ is given by
\begin{align*}
\Pi&=a+ar+ar^2+ar^3+\cdots\\
&=a\big(1+r+r^2+r^3+\cdots\big)\tag{2.1}
\end{align*}
Let $\lambda=1+r+r^2+r^3+\cdots$. In terms of $n$ for $n\geq0$, $\lambda$ can be expressed as
\begin{align*}
\lambda=\prod_{n=0}^{\infty}\big(1+r^{2^n}\big)=(1+r)\times(1+r^2)\times(1+r^4)\times(1+r^8)\times\cdots.\tag{2.2}
\end{align*}
Therefore, (2.1) follows that
\begin{align*}
\Pi=a(1+r)\times(1+r^2)\times(1+r^4)\times(1+r^8)\times\cdots\tag{2.3}
\end{align*}
We multiply the right hand side of (2.3) by $\displaystyle\frac{1-r}{1-r}$ to arrive at
\newpage
\begin{align*}
\Pi&=a\frac{\big[1-r\big]\big[1+r\big]\big[1+r^2\big]\big[1+r^4\big]\big[1+r^8\big]\big[\cdots\big]}{\big[1-r\big]}\\
&=a\bigg[\frac{1-r^{1+1+2+4+8+16+\cdots}}{1-r}\bigg]\\
\Pi&=a\bigg[\frac{r^{1+1+2+2^2+2^3+2^4+\cdots}-1}{r-1}\bigg]\tag{2.4}
\end{align*}
The expression in (2.4) in its current state does not look meaningful, but we can take the following approach. Let
\begin{align*}
\phi=1+2+2^2+2^3+\cdots\tag{2.5}
\end{align*}
so that 
\begin{align*}
\delta=1+\phi\tag{2.6}
\end{align*}
In terms of (2.6), we rewrite (2.4) as
\begin{align*}
\Pi=a\bigg(\frac{r^{\delta}-1}{r-1}\bigg)\tag{2.7}
\end{align*}
With reference to \cite{Kline} for instance, using the Euler's approach to divergent series one could obtain that the series in (2.5) is summable to -1, thus $\phi=-1$ and from (2.6) one sees that $\delta=0$. When one substitutes $\delta=0$ into (2.7) one obtains that $\Pi=0$, which means that the infinite geometric series in (2.5) is either not summable to -1 at all or is not  exactly summable to -1 as Euler had demonstrated. We then take the following approach: we have just shown that all infinite geometric series with the exception of the series with $r=1$ can be expressed of the form in (2.7). The infinite series in (2.5) is geometric with $a=1$ and $r=2$. On that note, substituting the above parameters into (2.7) gives
\begin{align*}
\Pi=\phi=a\bigg(\frac{r^{\delta}-1}{r-1}\bigg)=1\bigg(\frac{2^{\delta}-1}{2-1}\bigg)=2^{\delta}-1
\end{align*}
\begin{align*}
\phi+1=2^{\delta}\tag{2.8}
\end{align*}
But from (2.6), $\delta=\phi+1$. Therefore, (2.8) becomes
\begin{align*}
\phi+1=2^{\phi+1}\\
\implies& (\phi+1)e^{-\log2(\phi+1)}=1\\
\implies& -\log2(\phi+1)e^{-\log2(\phi+1)}=-\log2\\
\implies&-\log2(\phi+1)=W_n(-\log2)\\
\implies& \phi+1=-\frac{W_n(-\log2)}{\log2}.
\end{align*}
Therefore,
\begin{align*}
\delta=\phi+1=-\frac{W_n(-\log2)}{\log2}\tag{2.9}
\end{align*}
from which we see that 
\begin{align*}
\phi=\sum_{n=0}^{\infty}{2^n}=-1-\frac{W_n(-\log2)}{\log2}\tag{2.10}
\end{align*}
From (2.10) one sees that Euler's approach and the current approach are consistent in the sense that his -1 still forms part of (2.10). Euler offered a partial solution and now we give the exact sum of (2.5) in (2.10). Therefore, given that the constant ratio of an infinite geometric series is not equal to one, the exact sum of the series is determinable with
\begin{equation*}
\Pi=a\bigg(\frac{r^\delta-1}{r-1}\bigg),\,\,\, \delta=-\frac{W_n(-\log2)}{\log2}\tag{2.11}
\end{equation*}
as required.
\end{proof}
The following identities are important. From the transcendental equation
\begin{align*}
\phi+1=2^{\phi+1}\tag{2.12}
\end{align*}
we know from (2.9) that $\delta=\phi+1=-\frac{W_n(-\log2)}{\log2}$. Therefore,
\begin{align*}
\phi+1=2^{\phi+1}=-\frac{W_n(-\log2)}{\log2}=\delta=2^\delta\tag{2.13}
\end{align*}
From which one sees that
\begin{align*}
\log2=\frac{\log(1+\phi)}{\phi+1}=\frac{\log\delta}{\delta}=\frac{\log(1+1+2+2^2+2^3+\cdots)}{1+1+2+2^2+2^3+\cdots}\tag{2.14}
\end{align*}
Furthermore, from (2.12) we can rewrite
\begin{align*}
\phi+1&=2^{1+1+2+2^2+2^3+\cdots}\\
&=2^{1+\underbrace{1+2+2^2+2^3+\cdots}}\\
&=2^{1+1\big(\frac{2^{\delta}-1}{2-1}\big)}=2^{1+2^{\delta}-1}=2^{2^{\delta}}\\
&=2^{2^{1+\underbrace{1+2+2^2+\cdots}}}\\
&=2^{2^{2^{\delta}}}
\end{align*}
Repeating the process infinitely many times gives
\begin{align*}
\phi+1=2^{2^{2^{2^{2^{...}}}}}\tag{2.15}
\end{align*}
which leads us back to the fact that $\phi+1=-\frac{W_n(-\log2)}{\log2}$. But the idea in (2.15) is to show that
\begin{align*}
1+\sum_{n=0}^{\infty}{2^n}=2^{2^{2^{2^{2^{...}}}}}\tag{2.16}
\end{align*}
\subsection{\Large Some  general results}
Consider the following infinite series, whose finite sum we call $\kappa$:
\begin{align*}
\kappa=1+x+x^2+x^3+\cdots\tag{2.1.1}
\end{align*}
We do not assume that $|x|<1$, but $x\neq1$. Therefore, using the infinite geometric series formula in (2.11), we see that
\begin{align*}
\kappa=\frac{x^\delta-1}{x-1}
\end{align*}
From which we see using the quotient rule that  if $u=x^\delta-1$ and $v=x-1$, then
\begin{align*}
\frac{d\kappa}{dx}&=\frac{v\frac{du}{dx}-u\frac{dv}{dx}}{v^2}=\frac{\delta x^{\delta}-\delta x^{\delta-1}-x^{\delta}+1}{(x-1)^2}
\end{align*}
Therefore,
\begin{align*}
\frac{d\kappa}{dx}=1+2x+3x^2+4x^3+\cdots=\frac{\delta x^{\delta}-\delta x^{\delta-1}-x^{\delta}+1}{(x-1)^2}\tag{2.1.2}
\end{align*}
In section 3 it will become clear that in (2.1.2) for instance, when $|x|<1$ then $\delta x^{\delta}\approx\delta x^{\delta-1}\approx x^{\delta}\approx0$ and in such cases, one arrives at
\begin{align*}
1+2x+3x^2+4x^3+\cdots=\frac{1}{(x-1)^2}\tag{2.1.3}
\end{align*}
which is one of the commonest expressions in the literature. From (2.1.2) one further sees that multiplying both sides by $x$ gives
\begin{align*}
x+2x^2+3x^3+4x^4+\cdots=\frac{\delta x^{\delta+1}-\delta x^{\delta}-x^{\delta+1}+x}{(x-1)^2}\tag{2.1.4}
\end{align*}
\par Now suppose we integrate both sides of (2.1.1) with respect to $x$ from 0 to $m$ we will have
\begin{align*}
\int_{0}^{m}{(1+x+x^2+x^3+\cdots)}dx=\int_{0}^{m}{\frac{x^\delta-1}{x-1}}dx\\
m+\frac{m^2}{2}+\frac{m^3}{3}+\frac{m^4}{4}+\cdots=\int_{0}^{m}{\frac{x^\delta-1}{x-1}}dx\tag{2.1.5}
\end{align*}
The Wolfram online software gives,
\begin{align*}
\int_{0}^{m}{\frac{x^\delta-1}{x-1}}dx=-\frac{m^{\delta+1}}{\delta+1}{_2F_1}(1, \delta+1; \delta+2; m)-\log(1-m)\tag{2.1.6}
\end{align*}
where ${_2F_1}$ is a Gauss' hypergeometric series, which with reference to \cite{pearson}, is defined as
\begin{align*}
{_2F_1}(a,b;c;x)=\sum_{k=0}^{\infty}{\frac{(a)_k(b)_k}{(c)_kk!}}x^k
\end{align*}
where $(a)_k=\frac{\Gamma(a+k)}{\Gamma(a)}$, $\Gamma(a)=(a-1)!$, and $c\neq0, -1, -2,\cdots$.  On that note, we have
\begin{align*}
{_2F_1}(a,b;c;x)=\sum_{k=0}^{\infty}{\frac{\Gamma(a+k)\Gamma(b+k)\Gamma(c)}{\Gamma(a)\Gamma(b)\Gamma(c+k)}}\frac{x^k}{\Gamma(k+1)}\tag{2.1.7}
\end{align*}
By comparing (2.1.6) with (2.1.7), we see that  $x=m$, $a=1, b=\delta+1, c=\delta+2$, and substituting these parameters into (2.1.7) gives
\begin{align*}
-\frac{m^{\delta+1}}{\delta+1}{_2F_1}(1, \delta+1; \delta+2; m)&=-\frac{m^{\delta+1}}{\delta+1}\sum_{k=0}^{\infty}{\frac{\Gamma(1+k)\Gamma(\delta+1+k)\Gamma(\delta+2)}{\Gamma(1)\Gamma(\delta+1)\Gamma(\delta+2+k)}}\frac{m^k}{\Gamma(k+1)}\\
&=-\frac{m^{\delta+1}}{\delta+1}\sum_{k=0}^{\infty}{\frac{\delta+1}{\delta+1+k}}m^k\\
&=-\sum_{k=0}^{\infty}{\frac{m^{\delta+1+k}}{\delta+1+k}}
\end{align*}
Therefore, from (2.1.5) and (2.1.6) we write
\begin{align*}
m+\frac{m^2}{2}+\frac{m^3}{3}+\frac{m^4}{4}+\cdots=-\sum_{k=0}^{\infty}{\frac{m^{\delta+1+k}}{\delta+1+k}}-\log(1-m)
\end{align*}
\begin{align*}
-\log(1-m)=\sum_{k=0}^{\infty}{\frac{m^{\delta+1+k}}{\delta+1+k}}+m+\frac{m^2}{2}+\frac{m^3}{3}+\frac{m^4}{4}+\cdots
\end{align*}
\begin{align*}
\pi i-\log(m-1)=\sum_{k=0}^{\infty}{\frac{m^{\delta+1+k}}{\delta+1+k}}+m+\frac{m^2}{2}+\frac{m^3}{3}+\frac{m^4}{4}+\cdots
\end{align*}
Therefore,
\begin{align*}
\log(m-1)=\pi i-\sum_{k=0}^{\infty}{\frac{m^{\delta+1+k}}{\delta+1+k}}-\sum_{n=1}^{\infty}{\frac{m^n}{n}}\tag{2.1.8}
\end{align*}
From (2.1.8) when $|m|<1$, one sees that $\sum_{k=0}^{\infty}{\frac{m^{\delta+1+k}}{\delta+1+k}}\approx0$, simplifying (2.1.8) to
\begin{align*}
\log(1-m)=-m-\frac{m^2}{2}-\frac{m^3}{3}-\cdots\tag{2.1.9}
\end{align*}
which is already known in the literature. When $|m|>1$, then (2.1.8) would be defined in terms of the analytic continuation of the right hand side with respect to $m$. Furthermore, from
\begin{align*}
1+x+x^2+x^3+\cdots=\frac{x^\delta-1}{x-1}
\end{align*}
by the polynomial long division approach we have
\begin{align*}
1+x+x^2+x^3+\cdots=x^{\delta-1}+x^{\delta-2}+x^{\delta-3}+\cdots\tag{2.1.20}
\end{align*}
When $x\neq0$ in (2.1.20), integrating both sides with respect to $x$ from 0 to $m$ gives
\begin{align*}
m+\frac{m^2}{2}+\frac{m^3}{3}+\cdots=\frac{m^{\delta}}{\delta}+\frac{m^{\delta-1}}{\delta-1}+\frac{m^{\delta-2}}{\delta-2}+\cdots
\end{align*}
which follows that
\begin{align*}
1+\frac{m}{2}+\frac{m^2}{3}+\cdots=\frac{m^{\delta-1}}{\delta}+\frac{m^{\delta-2}}{\delta-1}+\frac{m^{\delta-3}}{\delta-2}+\cdots\tag{2.1.21}
\end{align*}
Similarly, differentiating both sides of (2.1.20) with respect to $x$ for all $x$ except $x=0$ gives
\begin{align*}
1+2x+3x^2+4x^3+\cdots&=x^{\delta-2}(\delta-1)+x^{\delta-3}(\delta-2)+x^{\delta-4}(\delta-3)+\cdots\\
&=\delta(x^{\delta-2}+x^{\delta-3}+x^{\delta-4}+\cdots)-(x^{\delta-2}+2x^{\delta-3}+3x^{\delta-4}+\cdots)\tag{2.1.22}
\end{align*}
But from (2.1.20), 
\begin{align*}
1+x+x^2+x^3+\cdots=x^{\delta-1}+x^{\delta-2}+x^{\delta-3}+\cdots
\end{align*}
Let $\tau(x)=1+x+x^2+x^3+\cdots$. Then
\begin{align*}
\tau(x)-x^{\delta-1}=x^{\delta-2}+x^{\delta-3}+x^{\delta-4}+\cdots\tag{2.1.23}
\end{align*}
Based on the result of (2.1.23), (2.1.22) becomes
\begin{align*}
\tau'(x)=\delta(\tau(x)-x^{\delta-1})-(x^{\delta-2}+2x^{\delta-3}+3x^{\delta-4}+\cdots)
\end{align*}
But since $\displaystyle\tau(x)=1+x+x^2+\cdots=\frac{x^\delta-1}{x-1}$, it follows from the previous step that
\begin{align*}
\tau'(x)&=\delta(\frac{x^{\delta}-1}{x-1}-x^{\delta-1})-(x^{\delta-2}+2x^{\delta-3}+3x^{\delta-4}+\cdots)\\
\tau'(x)&=\frac{x^{\delta}-x^{\delta}-1+x^{\delta-1}}{x-1}\delta-(x^{\delta-2}+2x^{\delta-3}+3x^{\delta-4}+\cdots)\\
\tau'(x)&=\frac{x^{\delta-1}-1}{x-1}\delta-(x^{\delta-2}+2x^{\delta-3}+3x^{\delta-4}+\cdots)
\end{align*}
But according to (2.1.2),
\begin{align*}
\tau'(x)=1+2x+3x^2+4x^3+\cdots=\frac{\delta x^{\delta}-\delta x^{\delta-1}-x^{\delta}+1}{(x-1)^2}
\end{align*}
Therefore,
\begin{align*}
\frac{\delta x^{\delta}-\delta x^{\delta-1}-x^{\delta}+1}{(x-1)^2}=\frac{x^{\delta-1}-1}{x-1}\delta-\sum_{n=2}^{\infty}{(n-1)x^{\delta-n}}
\end{align*}
\begin{align*}
\sum_{n=2}^{\infty}{(n-1)x^{\delta-n}}&=\frac{x^{\delta-1}-1}{x-1}\delta-\frac{\delta x^{\delta}-\delta x^{\delta-1}-x^{\delta}+1}{(x-1)^2}\\
&=\frac{\delta(x-1)(x^{\delta-1}-1)-(\delta x^{\delta}-\delta x^{\delta-1}-x^{\delta}+1)}{(x-1)^2}\\
&=\frac{\delta x^{\delta}-\delta x-\delta x^{\delta-1}+\delta-\delta x^{\delta}+\delta x^{\delta-1}+x^{\delta}-1}{(x-1)^2}\\
&=\frac{\delta-\delta x+x^{\delta}-1}{(x-1)^2},\,\,\, x\neq0,1
\end{align*}
Therefore, given that $x\neq0, \,\,1$, we have
\begin{align*}
x^{\delta-2}+2x^{\delta-3}+3x^{\delta-4}+\cdots=\frac{\delta-\delta x+x^{\delta}-1}{(x-1)^2}\tag{2.1.24}
\end{align*}
from which follows the result
\begin{align*}
1+\frac{2}{x}+\frac{3}{x^2}+\frac{4}{x^3}+\cdots=\frac{\delta-\delta x+x^{\delta}-1}{x^{\delta-2}(x-1)^2}, \,\, x\neq0, 1\tag{2.1.25}
\end{align*}
From (2.1.20), we equally remember that
\begin{align*}
1+x+x^2+x^3+\cdots&=x^{\delta-1}+x^{\delta-2}+x^{\delta-3}+\cdots\\
&=x^{\delta-1}\big(1+x^{-1}+x^{-2}+x^{-3}+x^{-4}+\cdots\big)
\end{align*}
Therefore,
\begin{align*}
\sum_{n=0}^{\infty}{x^n}=x^{\delta-1}\sum_{n=0}^{\infty}{x^{-n}}, \,\, x\neq0\tag{2.1.26}
\end{align*}
The relation in (2.1.26) then establishes the connection between an infinite geometric series and another infinite geometric series generated by taking the sum of the reciprocals of each of the terms of the original infinite geometric series. But we know from the geometric series formula in (2.11) that $\sum_{n=0}^{\infty}{x^n}=\frac{x^\delta-1}{x-1}$, so plugging that back into (2.1.26) and doing some rearrangements yields the following result:
\begin{align*}
\sum_{n=0}^{\infty}{x^{-n}}=\frac{\frac{x^\delta-1}{x-1}}{x^{\delta-1}}=\frac{x}{x-1}-\frac{x^{1-\delta}}{x-1}
\end{align*}
On that note we have
\begin{align*}
\sum_{n=0}^{\infty}{x^{-n}}=\frac{x}{x-1}-\frac{x^{1-\delta}}{x-1}\tag{2.1.27}
\end{align*}
We already know from the literature that in cases where an infinite geometric series converges, the limit of convergence becomes $\frac{a}{1-r}$. Therefore, in a case like the series in (2.1.27), using the formula $\frac{a}{1-r}$ we arrive at a limit of $\displaystyle\frac{x}{x-1}$ and comparing that to the result in (2.1.27) we see that using the formula $\frac{a}{1-r}$ introduces an error term $\displaystyle-\frac{x^{1-\delta}}{x-1}$. It turns out this error term is indeed very small and adding or subtracting it from the actual limit of the series does not really make a difference. Thus, if we rewrite (2.1.27) as 
\begin{align*}
1+\frac{1}{x}+\frac{1}{x^2}+\frac{1}{x^3}+\cdots=\frac{x}{x-1}\big(1-x^{-\delta}\big)
\end{align*}
and replace $\delta$ with $-\displaystyle\frac{W_n(-\log2)}{\log2}$, we obtain
\begin{align*}
1+\frac{1}{x}+\frac{1}{x^2}+\frac{1}{x^3}+\cdots=\frac{x}{x-1}\big(1-x^{--\frac{W_n(-\log2)}{\log2}}\big)=\frac{x}{x-1}\big(1-x^{\frac{W_n(-\log2)}{\log2}}\big)
\end{align*}
Therefore,
\begin{align*}
1+\frac{1}{x}+\frac{1}{x^2}+\frac{1}{x^3}+\cdots=\frac{x}{x-1}\big(1-x^{\frac{W_n(-\log2)}{\log2}}\big)\tag{2.1.28}
\end{align*}
Given that $x>0$, we see that as $n\rightarrow\pm\infty$, $\displaystyle x^{\frac{W_n(-\log2)}{\log2}}\approx0$. In cases where $x<0$ certain rules need to be followed in order to achieve the desired results and in section 3 we will learn more about such cases.\par
In the following section we will apply the above ideas to compute the sums of various infinite series. Most of these infinite series are very common in the literature. The idea is to reveal the consistency that exists between previous approaches to such infinite series and our present approach. What may be new, however, is our approach to some divergent series and not much is previously known about the sums of most of these divergent series.
\section{Some applications I: the sums of infinite series}
The point of departure in this section is the idea that for an infinite geometric series with a first term $a$ and a constant ratio $r$, its finite sum $\Pi$ is given by
\begin{align*}
\Pi=\big(r^\delta-1\big)\big(r-1\big)^{-1}, \,\, \delta=-W_n(-\log2)/\log2\tag{3.1}
\end{align*}
The nature of $W_n$ is such that at every $n$ there is a unique solution $\delta_n$ to the constant $\delta$. The solution $\delta_n$ gives a corresponding solution $\Pi_n$ as the finite sum of the series. Ideally, since the principal value of $W_n$ occurs at $n=0$, one could say that the principal values of $\delta_n$ and $\chi_n$ both occur at $n=0$. It turns out if we consider only the finite sum of the series at $n=0$, thus the solution $\Pi_0$, this solution alone is insufficient to entirely and accurately inform us about whether the series converges or diverges. But as we calculate the finite sums $\Pi_n$ at the other values of $n$, it becomes clear if a series converges and to what limit if it converges at all. \par The fact that at every $n$ there is a unique $\Pi_n$ implies that there are infinitely many $\Pi_{ns}$, all of which are complex in our case. These infinitely many complex solutions occur in conjugate pairs: one solution is $\Pi_n$ and its complex conjugate is $\Pi_{-(n+1)}$. If the series converges, all these infinitely many complex solutions will approach a single real solution as $n\rightarrow\pm\infty$. In fact, as $n\rightarrow\pm\infty$ the real parts of all these solutions will get closer and closer to exactly $\frac{a}{1-r}$ while their imaginary parts will approach zero at the same time. This consolidates the findings from past investigations that every convergent infinite geometric series has a finite sum of exactly $\frac{a}{1-r}$. In fact, previous investigations have already established that all infinite geometric series with $|r|<1$ converge while those with $|r|>1$ diverge and that is consistent with the findings of our current investigation. We find that the real parts of these  infinitely many complex solutions will approach $\frac{a}{1-r}$ and their imaginary parts will also approach zero as $n\rightarrow\pm\infty$  if and only if $|r|<1$. If however $|r|>1$, then as $n\rightarrow\pm\infty$ the real and the imaginary parts of these infinitely many complex solutions will keep growing bigger and bigger toward $\tilde\infty$. In cases where the series diverges, for practical purposes one may sidestep the reality for a moment and assume that the series' principal solution is the one at $n=0$. Cases involving infinite geometric series with $r<0$ are a bit more technical and require special attention and we will  explain more about such cases later. The sums of some non-geometric infinite series have  equally been discussed.\par
$\textbf{(I)}$. Let
\begin{align*}
\chi=1+\frac{1}{2}+\frac{1}{2^2}+\cdots\tag{3.2}.
\end{align*}
From the literature, the above series is known to converge to a finite sum of 2. We affirm that as follows. We apply the geometric series formula in (3.1) to arrive at 
\begin{align*}
\chi=a\bigg(\frac{r^{\delta}-1}{r-1}\bigg)=1\bigg(\frac{(2^{-1})^{\delta}-1}{2^{-1}-1}\bigg)=2\bigg(1-(2^{-1})^{\delta}\bigg)=2\bigg(1-\frac{1}{2^{\delta}}\bigg)
\end{align*}
Thus,
\begin{align*}
\chi=2\bigg(1-\frac{1}{2^{\delta}}\bigg)\tag{3.3}
\end{align*}
But we remember from (2.13) that $\delta=2^\delta=-\frac{W_n(-\log2)}{\log2}$. On that note (3.3) becomes
\begin{align*}
\chi=2\bigg(1-\frac{1}{2^{\delta}}\bigg)=2\bigg(1--\frac{\log2}{W_n(-\log2)}\bigg)=2+\frac{2\log2}{W_n(-\log2)}
\end{align*}
But at every $n$ there is a solution $\chi_n$. Therefore the n-th solution becomes
\begin{align*}
\chi_n=2+\frac{2\log2}{W_n(-\log2)}\tag{3.4}
\end{align*}
The following are some solutions of this series computed using the Wolfram online software.  Note how the real parts of these complex solutions get closer and closer to 2 while the imaginary parts approach 0 as $n\rightarrow\pm\infty$, and that is an indication that the series indeed converges to 2. Note also that these solutions come in conjugate pairs due to the nature of $W_n$.
\begin{multline*}
n=0; \chi_0=2+\frac{2\log2}{W_0(-\log2)}\\=1.4742143439565294465800856679560310297654568506549764091... \\-
0.99933887136663452730337627926333490785340383402837588188... i\\
n=-1; \chi_{-1}=2+\frac{2\log2}{W_{-1}(-\log2)}\\=1.4742143439565294465800856679560310297654568506549764091...\\ +0.99933887136663452730337627926333490785340383402837588188... i\\
n=1; \chi_1=2+\frac{2\log2}{W_1(-\log2)}\\=1.9462224290635914924225299071577906621194194158050879888...\\ -
0.16644763340153411531538593878155635981269719353730898578... i
\end{multline*}
\begin{multline*}
n=-2; \chi_{-2}=2+\frac{2\log2}{W_{-2}(-\log2)}\\=1.9462224290635914924225299071577906621194194158050879888...\\ +0.16644763340153411531538593878155635981269719353730898578... i\\
n=1000; \chi_{1000}=2+\frac{2\log2}{W_{1000}(-\log2)}\\=1.9999996801764241784705239085896264398069291735087686...\\ -
0.00022058004220891148156095076445122825171717796784197993... i\\
n=-1001; \chi_{-1001}=2+\frac{2\log2}{W_{-1001}(-\log2)}\\1.9999996801764241784705239085896264398069291735087686...\\ +0.00022058004220891148156095076445122825171717796784197993... i\\
n=100000; \chi_{100000}=2+\frac{2\log2}{W_{100000}(-\log2)}\\=1.999999999951831540897213403424847427024870806137407...\\ -
2.206350484675363065669178309195567447968599569711646...\times 10^{-6} i\\
n=-100001; \chi_{-100001}=2+\frac{2\log2}{W_{-100001}(-\log2)}\\=1.999999999951831540897213403424847427024870806137407...\\ +
2.206350484675363065669178309195567447968599569711646... \times 10^{-6} i\\
n=9999999999;\chi_{9999999999}=2+\frac{2\log2}{W_{9999999999}(-\log2)}\\=1.99999999999999999999114033873684937145991184975107...\\ -
2.20635600169199263374980310401207876144930029775140... \times10^{-11} i\\
n=-10000000000; \chi_{-10000000000}=2+\frac{2\log2}{W_{-10000000000}(-\log2)}\\=1.99999999999999999999114033873684937145991184975107...\\ +
2.20635600169199263374980310401207876144930029775140... \times 10^{-11} i
\end{multline*}
From the trend observed above, we conclude that indeed the infinite geometric series in (3.2) converges to 2. Since the series converges to 2 as $n\rightarrow\pm\infty$, we can write using a limit that
\begin{align*}
1+\frac{1}{2}+\frac{1}{2^2}+\cdots=\lim_{n\to\pm\infty}\bigg(2+\frac{2\log2}{W_n(-\log2)}\bigg)=2\tag{3.5}
\end{align*}
Thus, $\sum_{n=0}^{\infty}{2^{-n}}=2$.\par
$\textbf{(II)}$. Let $e$ be the euler's number with a value of $e=2.71828\cdots$. We investigate whether the following  infinite geometric series converges or diverges.
\begin{align*}
\Gamma=1+e+e^2+e^3+\cdots.\tag{3.6}
\end{align*}
The first term $a=1$ and the constant ratio $r=e$ and plugging these parameters into (3.1) gives
\begin{align*}
\Gamma=\frac{e^{-\frac{W_n(-\log2)}{\log2}}-1}{e-1}\tag{3.7}
\end{align*}
But since every $n$ gives gives a corresponding $\Gamma_n$ we write
\begin{align*}
\Gamma_n=\frac{e^{-\frac{W_n(-\log2)}{\log2}}-1}{e-1}\tag{3.8}
\end{align*}
We calculate some solutions of (3.8) below. This time around, note how the real and the imaginary components of the solutions keep growing bigger and bigger as $n\rightarrow\pm\infty$. In other words, the series in (3.6) does not converge. 
\begin{multline*}
n=0; \Gamma_0=\big(\frac{e^{-\frac{W_0(-\log2)}{\log2}}-1}{e-1}\big)\\=-0.5775104848114135400681245842452675099759474284944265796...\\ -
1.327567118044092590492574213638171838585129990371402738... i\\
n=-1; \Gamma_{-1}=\big(\frac{e^{-\frac{W_{-1}(-\log2)}{\log2}}-1}{e-1}\big)\\=-0.5775104848114135400681245842452675099759474284944265796...\\ +
1.327567118044092590492574213638171838585129990371402738... i\\
n=1; \Gamma_{1}=\big(\frac{e^{-\frac{W_1(-\log2)}{\log2}}-1}{e-1}\big)\\=-2.8375054648584783816269775165317195254290405213806312504... \\+
19.437893097403337265141872906251855541488893958220547286... i\\
n=-2; \Gamma_{-2}=\big(\frac{e^{-\frac{W_{-2}(-\log2)}{\log2}}-1}{e-1}\big)\\=-2.8375054648584783816269775165317195254290405213806312504...\\ -
19.437893097403337265141872906251855541488893958220547286... i\\
n=2; \Gamma_2=\big(\frac{e^{-\frac{W_2(-\log2)}{\log2}}-1}{e-1}\big)\\=14.3348668396616284138248078732154747631931821005012981975...\\ -
43.1007525458601449184213803693033978514259045444876992179... i
\end{multline*}
\begin{multline*}
n=-3; \Gamma_{-3}=\big(\frac{e^{-\frac{W_{-3}(-\log2)}{\log2}}-1}{e-1}\big)\\=14.3348668396616284138248078732154747631931821005012981975... \\+
43.1007525458601449184213803693033978514259045444876992179... i\\
n=1000; \Gamma_{1000}=280204.803597047442334074027807659813954790115677357221173...\\ -
101639.981827083908652521726827728473646144121522199885534... i\\
n=-1001; \Gamma_{-1001}=\big(\frac{e^{-\frac{W_{-1001}(-\log2)}{\log2}}-1}{e-1}\big)\\=280204.803597047442334074027807659813954790115677357221173...\\ +
101639.981827083908652521726827728473646144121522199885534... i\\
n=999999999; \Gamma_{999999999}=\\=4.727589932399994187919828134433908858696081877476299...\times 10^{13}\\ +
1.264507522528720682623483772861042857074729991130866... \times 10^{14 }\\
n=-1000000000; \Gamma_{-1000000000}\\=4.727589932399994187919828134433908858696081877476299... \times10^{13}\\ -
1.264507522528720682623483772861042857074729991130866... \times 10^{14 }i
\end{multline*}
As $n\rightarrow\pm\infty$, $\Gamma_n\rightarrow\tilde\infty$. Therefore, this is a proof that the series diverges.\par
$\textbf{(III)}$. In a similar way, we consider
\begin{align*}
\tau=1+\frac{1}{5}+\frac{1}{5^2}+\frac{1}{5^3}+\cdots.\tag{3.9}
\end{align*}
Substituting $a=1$ and $r=5^{-1}$ into the geometric series formula gives
\begin{align*}
\tau&=1\bigg(\frac{(5^{-1})^{\delta}-1}{5^{-1}-1}\bigg)=\frac{5}{4}\bigg(1-(5^{-1})^{\delta}\bigg)\\
\tau&=\frac{5}{4}\big(1-(5^{-1})^{-\frac{W_n(-\log2)}{\log2}}\big)\\
\tau_n&=\frac{5}{4}\big(1-5^{\frac{W_n(-\log2)}{\log2}}\big)\tag{3.10}
\end{align*}
Based on previous knowledge, we know that the geometric series in (3.9) converges to 1.25. Below are some of the solutions of the same series from our new perspective. Note how these solutions get closer and closer to 1.25 as $n\rightarrow\pm\infty$.
\begin{multline*}
n=0; \tau_0=\frac{5}{4}\bigg(1-5^{\frac{W_0(-\log2)}{\log2}}\bigg)\\=1.5200114549575887328456395565619943089723076226224848404...\\ -
0.19231881443558700658293287982219050748031523776629203782... i
\end{multline*}
\begin{multline*}
n=-1; \tau_{-1}=\frac{5}{4}\bigg(1-5^{\frac{W_{-1}(-\log2)}{\log2}}\bigg)\\=1.5200114549575887328456395565619943089723076226224848404...\\ +
0.19231881443558700658293287982219050748031523776629203782... i\\
n=1; \tau_1\\=1.24899661085255744238509077952971084688169565073734812...\\ +
0.00424687810789892661799570484371450277687000990178762897... i\\
n=-2; \tau_{-2}=\frac{5}{4}\big(1-5^{\frac{W_{-2}(-\log2)}{\log2}}\big)\\=1.24899661085255744238509077952971084688169565073734812...\\ -
0.00424687810789892661799570484371450277687000990178762897... i\\
n=1000; \tau_{1000}=\frac{5}{4}\big(1-5^{\frac{W_{1000}(-\log2)}{\log2}}\big)\\=1.25000000080799784295226199540861938591041232751721...\\ +
4.08578217219123188734549651669459598127260070766517... \times 10^{-11} i\\
n=-1001; \tau_{-1001}=\frac{5}{4}\big(1-5^{\frac{W_{-1001}(-\log2)}{\log2}}\big)\\=1.25000000080799784295226199540861938591041232751721...\\ -
4.08578217219123188734549651669459598127260070766517... \times 10^{-11} i\\
n=99999; \tau_{99999}=\frac{5}{4}\big(1-5^{\frac{W_{99999}(-\log2)}{\log2}}\big)\\=1.24999999999998327300321452684437229521927020607450...\\ -
7.62027159648786072814483552678385138721898360430759... \times 10^{-15}i \\
n=-100000; \tau_{-100000}=\frac{5}{4}\big(1-5^{\frac{W_n(-\log2)}{\log2}}\big)\\=1.24999999999998327300321452684437229521927020607450... \\+
7.62027159648786072814483552678385138721898360430759... \times 10^{-15} i\\
n=999999999; \tau_{999999999}=\frac{5}{4}\big(1-5^{\frac{W_{999999999}(-\log2)}{\log2}}\big)\\=1.24999999999999999999999423497657167048361241213229...\\ -
7.52123733856768524780038524459798608921127211841592... \times10^{-24} i\\
n=-1000000000;\tau_{-1000000000}=\frac{5}{4}\big(1-5^{\frac{W_{-1000000000}(-\log2)}{\log2}}\big)\\=1.24999999999999999999999423497657167048361241213229...\\ +
7.52123733856768524780038524459798608921127211841592... \times10^{-24} i
\end{multline*}
The few solutions above prove that indeed the series in (3.9) has a limit of convergence of 1.25. Therefore, we may write
\begin{align*}
\tau=1+\frac{1}{5}+\frac{1}{5^2}+\frac{1}{5^3}+\cdots=\lim_{n\to\pm\infty}\big(\frac{5}{4}\big(1-5^{\frac{W_n(-\log2)}{\log2}}\big)\big)=1.25\tag{3.11}
\end{align*}
\par So far all the infinite geometric series we have considered have involved those with constant ratios greater than zero, thus those with $r>0$. For infinite geometric series with negative constant ratios, $r<0$, things can get a bit more complicated and confusing if care is not taken. It is already established in the literature that any infinite geometric series of which the inequality $|r|>1$ is true diverges and converges for $|r|<1$. Note here that $r$ is the constant ratio of the series. So for an infinite geometric series with with $r=-2$, for instance, that series diverges since $|-2|=2>1$. Assuming the first term of that series is 1, then by definition we generate the series as\par
$\textbf{(IV)}$.
\begin{align*}
1-2+2^2-2^3+\cdots.\tag{3.12}
\end{align*}
Using approaches like the Euler summation (transform) and the formula $\frac{a}{1-r}$, one arrives at a  finite sum of $1/3$ for the series in (3.12), but in reality this series diverges. Let $A$ be the finite sum of (3.12). Putting $a=1$ and $r=-2$ into the geometric series formula in (3.1) gives
\begin{align*}
A_n&=1\big(\frac{(-2))^{\delta}-1}{-2-1}\big)=\frac{1}{3}\big(1-(-2)^{\delta}\big)\\
A_n&=\frac{1}{3}\big(1-(-2)^{-\frac{W_n(-\log2)}{\log2}}\big)\tag{3.13}
\end{align*}
Some solutions of (3.13) are shown below. We expected the imaginary and the real parts of all these solutions to grow bigger and bigger as $n\rightarrow\pm\infty$ since $|-2|>1$, but that is not what we see. Some solutions get closer and closer to $1/3$ as predicted by the Euler transform while others diverge. But all of them were supposed to diverge since $|-2|=2>1$.
\begin{multline*}
n=0; A_0=\frac{1}{3}\big(1-(-2)^{-\frac{W_0(-\log2)}{\log2}}\big)\\=-5.0604898851734126859407445402009994471444055594104308199...\\ -
81.047460084156704566079711546098232685489423179940352285... i\\
n=-1; A_{-1}=\frac{1}{3}\big(1-(-2)^{-\frac{W_{-1}(-\log2)}{\log2}}\big)\\=0.337023323787748744597344511022967822937307465616495830...\\ +
0.00219023979701119344466133805752647859087940694772316425... i\\
n=1; A_1=\frac{1}{3}\big(1-(-2)^{-\frac{W_1(-\log2)}{\log2}}\big)\\=2.493015956698306661097353014865712618641973709445937... \times10^{15}\\ +
9.394344172313893836256957030436378146782153697839025... \times10^{14} i
\end{multline*}
\begin{multline*}
n=-2; A_{-2}=\frac{1}{3}\big(1-(-2)^{-\frac{W_{-2}(-\log2)}{\log2}}\big)\\=0.333333333333328070798996220170277053010943402202120...\\ +
1.42610275184498257407923062410667357077524582026274... \times10^{-15} i\\
n=100; A_{100}=\frac{1}{3}\big(1-(-2)^{-\frac{W_{100}(-\log2)}{\log2}}\big)\\=1.05898855121434388418793778650276144287983290623472... \times10^{1242}\\ +
1.80645545000327589428224761159850219061969044259878...\times10^{1242} i\\
n=-101; A_{-101}=\frac{1}{3}\big(1-(-2)^{-\frac{W_{-101}(-\log2)}{\log2}}\big)\\=0.33333333333333333333333333333333333333333333333333...\\ -
3.7317303553187626438395160129112698504347014759232...\times10^{-1238} i
\end{multline*}
We consider a second case with $r<0$.\par
$\textbf{(V)}$. We consider the alternating geometric series 
\begin{align*}
1-\frac{1}{2}+\frac{1}{2^2}-\frac{1}{2^3}+\cdots.\tag{3.14}
\end{align*}
whose finite sum we denote by $\mu$. Plugging the parameters $a=1$ and $r=-2^{-1}$ into the geometric series formula gives
\begin{align*}
\mu_n&=\frac{2}{3}\big(1-(-2^{-1})^{\delta}\big)\\
\mu_n&=\frac{2}{3}\big(1-(-2^{-1})^{-\frac{W_n(-\log2)}{\log2}}\big)\\
\mu_n&=\frac{2}{3}\big(1-(-2)^{\frac{W_n(-\log2)}{\log2}}\big)\tag{3.15}
\end{align*}
The following are few solutions of the series.
\begin{multline*}
\mu_0=\frac{2}{3}\big(1-(-2)^{\frac{W_0(-\log2)}{\log2}}\big)\\=0.666484995462972031527452010960193252705973807426014424...\\ +
0.00272978720907317928266677011662587036421982364032329564... i\\
\mu_{-1}=\frac{2}{3}\big(1-(-2)^{\frac{W_{-1}(-\log2)}{\log2}}\big)\\=45.1998696347921884356566880183651911570171387276375611093...\\ -
26.4332373305876459271297539349528138148246004008534579496... i\\
\mu_1=\frac{2}{3}\big(1-(-2)^{\frac{W_1(-\log2)}{\log2}}\big)\\=0.666666666666666744720998661760143468299145605078681... \\-
2.94129388515060068436890136878793587165027519899709... \times10^{-17} i
\end{multline*}
\begin{multline*}
\mu_{-2}=\frac{2}{3}\big(1-(-2)^{\frac{W_{-2}(-\log2)}{\log2}}\big)\\=-3.93383554520029633859385310565266226635455157791543... \times10^{13}\\ -
1.06603650198571788707311255266164424367635506502647...\times 10^{13} i\\
\mu_{20000}=0.6666666666666666666666666666666666666666666666666...\\ +
4.729886948766748335989907468272251834228378563709...\times 10^{-247364} i\\
\mu_{-20001}=\frac{2}{3}\big(1-(-2)^{\frac{W_{-20001}(-\log2)}{\log2}}\big)\\=-2.852544506329290711818524569134494370266629729544...\times 10^{247351}\\ -
2.880807750869019223752711736357069580065465999786...\times 10^{247350} i
\end{multline*}
Undoubtedly, all the solutions of (3.14) were supposed to converge to $2/3$ since $|-2^{-1}|<1$. Why could that be happening? It turns out for $r<0$, the problem originates from the term $r^\delta$.  From the geometric series formula, we remember that
\begin{align*}
\Pi=a\big[r^\delta-1\big]\big[r-1\big]^{-1}, \,\,\delta=-W_n(-\log2)/\log2
\end{align*}
We equally know that all the possible solutions of $\delta$ are complex. Let $r=-x$. Therefore
\begin{align*}
\Pi=a\big[(-x)^\delta-1\big]\big[-x-1\big]^{-1}=\frac{a(1-(-x)^\delta)}{x+1}\tag{3.16}
\end{align*}
We remember that the infinite series for $\delta$ is
\begin{align*}
\delta=1+1+2+2^2+2^3+\cdots=-\frac{W_n(-\log2)}{\log2}
\end{align*}
So in a case like $(-x)^\delta$ as in (3.16), what it means is that
\begin{align*}
(-x)^\delta=(-x)^{1+1+2+2^2+2^3+\cdots}=x^{1+1+2+2^2+2^3+\cdots}
\end{align*}
Since $(-x)^{1+1+2+2^2+2^3+\cdots}=x^{1+1+2+2^2+2^3+\cdots}$ then $(-x)^\delta$ should also be treated as $x^\delta$. It turns out if we follow this convention, then all the anomalies encountered in the last two cases would be rectified. Therefore, from (3.16) we now have $\Pi=\frac{a(1-(-x)^\delta)}{x+1}=\frac{a(1-x^\delta)}{x+1}$. \par
We revisit the case in (IV). We recall that the formula for the sum of the series in (3.12) as we saw from (3.13) was $A_n=\frac{1}{3}\big(1-(-2)^{-\frac{W_n(-\log2)}{\log2}}\big)$, which we now rewrite  as
\begin{align*}
A_n=\frac{1}{3}\big(1-2^{-\frac{W_n(-\log2)}{\log2}}\big)\tag{3.17}
\end{align*}
We compute some of the solutions of (3.17) below. Remember, at first some of the solutions converged while others diverged, and those solutions never occurred in conjugate pairs. Note how all the solutions now diverge and exist in conjugate pairs as usual.
\begin{multline*}
A_0=\frac{1}{3}\big(1-2^{-\frac{W_0(-\log2)}{\log2}}\big)\\=0.05844048461930859256197846853946558667513085993945639110...\\ +
0.5224773746165492870195247970643289758444358580680725210... i\\
A_{-1}=\frac{1}{3}\big(1-2^{-\frac{W_{-1}(-\log2)}{\log2}}\big)\\=0.05844048461930859256197846853946558667513085993945639110...\\ -
0.5224773746165492870195247970643289758444358580680725210... i\\
A_1=\frac{1}{3}\big(1-2^{-\frac{W_1(-\log2)}{\log2}}\big)\\=-0.8384122406388703761769727134604899532429586726628541870... \\+
3.626684402815846989110558653650412782933100110335410222... i\\
A_{-2}=-0.8384122406388703761769727134604899532429586726628541870...\\ -
3.626684402815846989110558653650412782933100110335410222... i\\
A_{20000}=\frac{1}{3}\big(1-2^{-\frac{W_{20000}(-\log2)}{\log2}}\big)\\=-5.4893309737107893178531458699866562459241155714946987482... \\+
60432.223904718177630874199247065308240623466813306269892... i\\
A_{-20001}=\frac{1}{3}\big(1-2^{-\frac{W_{-20001}(-\log2)}{\log2}}\big)\\=-5.4893309737107893178531458699866562459241155714946987482... \\-
60432.223904718177630874199247065308240623466813306269892... i\\
A_{999999999}=\frac{1}{3}\big(1-2^{-\frac{W_{999999999}(-\log2)}{\log2}}\big)\\=-10.69253845177381875693948493479140877814806147798654...\\ +
3.021573425618615800416703061183890012675764968419242...\times 10^9 i\\
A_{-1000000000}=\frac{1}{3}\big(1-2^{-\frac{W_{-1000000000}(-\log2)}{\log2}}\big)\\=-10.69253845177381875693948493479140877814806147798654...\\ -
3.021573425618615800416703061183890012675764968419242... \times 10^9 i
\end{multline*}
Also, we revist the case in (V) and using the same approach we will find that the formula for the sum of the series in (3.14) as seen in (3.15) now becomes
\begin{align*}
\mu_n=\frac{2}{3}\big(1-2^{\frac{W_n(-\log2)}{\log2}}\big)\tag{3.18}
\end{align*}
Some solutions of (3.18) are
\begin{multline*}
\mu_0=\frac{2}{3}(1-2^{\frac{W_0(-\log2)}{\log2}})\\=0.49140478131884314886002855598534367658848561688499213639...\\ -0.33311295712221150910112542642111163595113461134279196062... i\\
\mu_{-1}=\frac{2}{3}(1-2^{\frac{W_{-1}(-\log2)}{\log2}})\\=0.49140478131884314886002855598534367658848561688499213639...\\ +
0.33311295712221150910112542642111163595113461134279196062... i\\
\mu_1=\frac{2}{3}(1-2^{\frac{W_1(-\log2)}{\log2}})\\=0.6487408096878638308075099690525968873731398052683626629...\\ -
0.05548254446717803843846197959385211993756573117910299526... i\\
\mu_{-2}=\frac{2}{3}(1-2^{\frac{W_{-2}(-\log2)}{\log2}})\\=0.6487408096878638308075099690525968873731398052683626629...\\ +
0.05548254446717803843846197959385211993756573117910299526... i\\
\mu_{20000}=\frac{2}{3}(1-2^{\frac{W_{20000}(-\log2)}{\log2}})\\=0.6666666663123659085784802828183925456718609879460362...\\ -
3.677214006051134174398260962361881494702379797567674... \times 10^{-6} i\\
\mu_{-20001}=\frac{2}{3}(1-2^{\frac{W_{-20001}(-\log2)}{\log2}})\\=0.6666666663123659085784802828183925456718609879460362...\\ +
3.677214006051134174398260962361881494702379797567674...\times 10^{-6} i\\
\mu_{999999999}=\frac{2}{3}(1-2^{\frac{W_{999999999}(-\log2)}{\log2}})\\=0.666666666666666666398296571833290755514913440924205...\\ -
7.35452001060427636084341477053333302747723398650694... \times10^{-11} i\\
\mu_{-1000000000}=\frac{2}{3}(1-2^{\frac{W_{-1000000000}(-\log2)}{\log2}})\\=0.666666666666666666398296571833290755514913440924205... \\+
7.35452001060427636084341477053333302747723398650694... \times10^{-11} i
\end{multline*}
All the solutions of this series too now converge and occur in conjugate pairs as expected. As $n\rightarrow\pm\infty$, $\mu_n\rightarrow2/3$.  In that case we write

\begin{align*}
\mu=\lim_{n\to\pm\infty}(\frac{2}{3}(1-2^{\frac{W_n(-\log2)}{\log2}}))=2/3\tag{3.19}
\end{align*}
$\textbf{(VI)}$. Consider the following infinite series:
\begin{align*}
\sum_{n=1}^{\infty}{n2^{n-1}}=1+4+12+32+\cdots.\tag{3.20}
\end{align*}
This is a slightly different case from those we have previously dealt with, since the current infinite series is non-geometric. We remember the following relation from (2.1.2).
\begin{align*}
1+2x+3x^2+4x^3+\cdots=\frac{\delta x^{\delta}-\delta x^{\delta-1}-x^{\delta}+1}{(x-1)^2}\tag{3.21}
\end{align*}
By comparing (3.20) with (3.21) we see that $x=2$. Therefore,
\begin{align*}
1+4+12+32+\cdots&=\frac{\delta 2^{\delta}-\delta 2^{\delta-1}-2^{\delta}+1}{(2-1)^2}\\
&=\delta 2^{\delta}-\delta 2^{\delta-1}-2^{\delta}+1\\
&=\delta2^{\delta-1}-2^{\delta}+1\\
&=\frac{\delta2^{\delta}}{2}-2^{\delta}+1\tag{3.22}
\end{align*}
But we recall from (2.13) that $\delta=2^\delta=-\frac{W_n(-\log2)}{\log2}$ and subsituting that into (3.22) gives
\begin{align*}
1+4+12+32+\cdots&=\frac{\delta2^{\delta}}{2}-2^{\delta}+1\\
&=\frac{1}{2}\bigg[-\frac{W_n(-\log2)}{\log2}\bigg]^2-\bigg[-\frac{W_n(-\log2)}{\log2}\bigg]+1\\
&=\frac{1}{2}\bigg[\frac{W_n(-\log2)}{\log2}\bigg]^2+\bigg[\frac{W_n(-\log2)}{\log2}\bigg]+1
\end{align*}
Therefore,
\begin{align*}
\sum_{n=1}^{\infty}{n2^{n-1}}=\frac{1}{2}\bigg[\frac{W_n(-\log2)}{\log2}\bigg]^2+\bigg[\frac{W_n(-\log2)}{\log2}\bigg]+1\tag{3.23}
\end{align*}
If we were to compute some few solutions of this series as usual, we would find that as $n\rightarrow\pm\infty$ the real and the imaginary parts of the solutions would keep growing bigger and bigger toward $\tilde\infty$. This is an indication that the series diverges. \par
$\textbf{(VII)}$.  Let $\varphi$ denote the finite sum of the infinite series
\begin{align*}
2+\frac{5}{2}+\frac{17}{4}+\frac{65}{8}+\frac{257}{16}+\frac{1025}{32}+\cdots\tag{3.24}
\end{align*}
which can be expressed as
\begin{align*}
\varphi=\sum_{n=0}^{\infty}\bigg({\frac{2^{2n}+1}{2^n}}\bigg)=\sum_{n=0}^{\infty}{2^n}+\sum_{n=0}^{\infty}{2^{-n}}.\tag{3.25}
\end{align*}
We have already shown that
\begin{align*}
\phi=\sum_{n=0}^{\infty}{2^n}=-1-\frac{W_n(-\log2)}{\log2}
\end{align*}
and that
\begin{align*}
\sum_{n=0}^{\infty}{2^{-n}}=2+\frac{2\log2}{W_n(-\log2)}
\end{align*}
Adding these two results gives
\begin{align*}
\varphi=\sum_{n=0}^{\infty}{2^n}+\sum_{n=0}^{\infty}{2^{-n}}=-1-\frac{W_n(-\log2)}{\log2}+2+\frac{2\log2}{W_n(-\log2)}
\end{align*}
\begin{align*}
\varphi=1+\frac{2\log2}{W_n(-\log2)}-\frac{W_n(-\log2)}{\log2}\tag{3.26}
\end{align*}
Some solutions of the series are:
\begin{multline*}
n=0; \varphi_0=1+\frac{2\log2}{W_0(-\log2)}-\frac{W_0(-\log2)}{\log2}\\=1.29889289009860366889415026233763426974006427083660723587...\\ -
2.56677099521628238836195067045632183538671140823259344505... i\\
n=-1; \varphi_{-1}=1+\frac{2\log2}{W_{-1}(-\log2)}-\frac{W_{-1}(-\log2)}{\log2}\\=1.29889289009860366889415026233763426974006427083660723587...\\ +
2.56677099521628238836195067045632183538671140823259344505... i\\
n=1; \varphi_1=1+\frac{2\log2}{W_1(-\log2)}-\frac{W_1(-\log2)}{\log2}\\=4.46145915098020262095344804753926052184829543379365055008...\\ -
11.0465008418490750826470618997327947086119975245435396532... i
\end{multline*}
\begin{multline*}
n=-2; \varphi_{-2}=1+\frac{2\log2}{W_{-2}(-\log2)}-\frac{W_{-2}(-\log2)}{\log2}\\=4.46145915098020262095344804753926052184829543379365055008...\\ +
11.0465008418490750826470618997327947086119975245435396532... i\\
n=1000; \varphi_{1000}=1+\frac{2\log2}{W_{1000}(-\log2)}-\frac{W_{1000}(-\log2)}{\log2}\\=14.1464082795691260528587498535118584393446550864750150266...\\ -
9066.98459251336241786543537539012639711180115904678877743... i\\
n=-1001; \varphi_{-1001}=1+\frac{2\log2}{W_{-1001}(-\log2)}-\frac{W_{-1001}(-\log2)}{\log2}\\=14.1464082795691260528587498535118584393446550864750150266...\\ +
9066.98459251336241786543537539012639711180115904678877743... i\\
n=999999999; \varphi_{999999999}=1+\frac{2\log2}{W_{999999999}(-\log2)}-\frac{W_{999999999}(-\log2)}{\log2}\\=34.077615355321456270013344519874098600988924756732256...\\ -
9.0647202768558474014707447838697983288525973483737266... \times10^9 i\\
n=-1000000000; \varphi_{-1000000000}=1+\frac{2\log2}{W_{-1000000000}(-\log2)}-\frac{W_{-1000000000}(-\log2)}{\log2}\\=34.077615355321456270013344519874098600988924756732256...\\ +
9.0647202768558474014707447838697983288525973483737266... \times10^9 i
\end{multline*}
From these solutions we can say that since the real and the imaginary parts of the solutions grow bigger toward $\tilde\infty$  as $n\rightarrow\pm\infty$, the series diverges.\par
$\textbf{(VIII)}$. We consider
\begin{align*}
k=\frac{1}{8}+\frac{1}{8^2}+\frac{1}{8^3}+\frac{1}{8^4}+\cdots\tag{3.27}.
\end{align*}
\begin{align*}
k&=a\big(\frac{r^{\delta}-1}{r-1}\big)=\frac{1}{8}\big(\frac{(8^{-1})^{\delta}-1}{8^{-1}-1}\big)\\
&=\frac{1}{7}\big(1-(8^{-1})^{\delta}\big)\\
k_n&=\frac{1}{7}\big(1-(8^{-1})^{-\frac{W_n(-\log2)}{\log2}}\big)\\
k_n&=\frac{1}{7}\big(1-8^{\frac{W_n(-\log2)}{\log2}}\big)\tag{3.28}
\end{align*}
Some solutions of (3.28) are:

\begin{multline*}
k_0=\frac{1}{7}(1-8^{\frac{W_0(-\log2)}{\log2}})\\=0.168391397301116745552441557679736398538306784464523582... \\+
0.00302168865921421821896599738863331053309698223034262169... i\\
k_{-1}=\frac{1}{7}(1-8^{\frac{W_{-1}(-\log2)}{\log2}})\\=0.168391397301116745552441557679736398538306784464523582...\\ -0.00302168865921421821896599738863331053309698223034262169... i\\
k_{10000000}=\frac{1}{7}(1-8^{\frac{W_{10000000}(-\log2)}{\log2}})\\=0.142857142857142857142857142857310646523718725897358...\\ +
1.91795629021040058978005704291971341598096413468966...\times10^{-25} i\\
k_{-10000001}=\frac{1}{7}(1-8^{\frac{W_{-10000001}(-\log2)}{\log2}})\\=0.142857142857142857142857142857310646523718725897358...\\ -
1.91795629021040058978005704291971341598096413468966...\times 10^{-25} i
\end{multline*}
As $n\rightarrow\pm\infty$, $k_n\rightarrow0.\overline{142857}=1/7$. Therefore,
\begin{align*}
k=\lim_{n\to\pm\infty}(\frac{1}{7}(1-8^{\frac{W_n(-\log2)}{\log2}}))=\frac{1}{7}\tag{3.29}
\end{align*}
$\textbf{(IX)}$. Suppose we were to evaluate the following:
\begin{align*}
\rho=\frac{\sqrt{2}+2(\sqrt{2})^2+3(\sqrt{2})^3+4(\sqrt{2})^4+\cdots}{1+\sqrt{2}+(\sqrt{2})^2+(\sqrt{2})^3+(\sqrt{2})^4+\cdots}
\end{align*}
what would $\rho$ be? We remember from (2.1.4) that
\begin{align*}
x+2x^2+3x^3+4x^4+\cdots=\frac{\delta x^{\delta+1}-\delta x^{\delta}-x^{\delta+1}+x}{(x-1)^2}
\end{align*}
So when $x=\sqrt{2}$, we have
\begin{align*}
\sqrt{2}+2(\sqrt{2})^2+3(\sqrt{2})^3+4(\sqrt{2})^4+\cdots=\frac{\delta (\sqrt{2})^{\delta+1}-\delta (\sqrt{2})^{\delta}-(\sqrt{2})^{\delta+1}+\sqrt{2}}{(\sqrt{2}-1)^2}\tag{3.30}
\end{align*}
We equally remember from the geometric series formula that
\begin{align*}
1+x+x^2+x^3+\cdots=\frac{x^\delta-1}{x-1}
\end{align*}
So when $x=\sqrt{2}$, we have
\begin{align*}
1+\sqrt{2}+(\sqrt{2})^2+(\sqrt{2})^3+(\sqrt{2})^4+\cdots=\frac{(\sqrt{2})^\delta-1}{\sqrt{2}-1}\tag{3.31}
\end{align*}
Dividing (3.30) by (3.31) gives
\begin{align*}
\rho=\frac{\sqrt{2}+2(\sqrt{2})^2+3(\sqrt{2})^3+4(\sqrt{2})^4+\cdots}{1+\sqrt{2}+(\sqrt{2})^2+(\sqrt{2})^3+(\sqrt{2})^4+\cdots}=\frac{\frac{\delta (\sqrt{2})^{\delta+1}-\delta (\sqrt{2})^{\delta}-(\sqrt{2})^{\delta+1}+\sqrt{2}}{(\sqrt{2}-1)^2}}{\frac{(\sqrt{2})^\delta-1}{\sqrt{2}-1}}
\end{align*}
Therefore,
\begin{align*}
\rho=\frac{\sqrt{2}+2(\sqrt{2})^2+3(\sqrt{2})^3+4(\sqrt{2})^4+\cdots}{1+\sqrt{2}+(\sqrt{2})^2+(\sqrt{2})^3+(\sqrt{2})^4+\cdots}&=\frac{\delta(\sqrt{2})^{\delta+1}-\delta (\sqrt{2})^{\delta}-(\sqrt{2})^{\delta+1}+\sqrt{2}}{(\sqrt{2}-1)((\sqrt{2})^\delta-1)}\\
&=\frac{(1+\sqrt{2})\big(\delta\sqrt{2}^{\delta+1}-\delta \sqrt{2}^{\delta}-\sqrt{2}^{\delta+1}+\sqrt{2}\big)}{\sqrt{2}^\delta-1}
\end{align*}
But we remember from (2.13) that $\delta=2^\delta=-\frac{W_n(-\log2)}{\log2}$. On that note, we obtain by simplifying the last step above using the basic laws of indices that
\begin{align*}
\rho=\frac{\sqrt{2}+2(\sqrt{2})^2+3(\sqrt{2})^3+4(\sqrt{2})^4+\cdots}{1+\sqrt{2}+(\sqrt{2})^2+(\sqrt{2})^3+(\sqrt{2})^4+\cdots}=\frac{(1+\sqrt{2})\big(\delta\sqrt{2\delta}-\delta\sqrt{\delta}-\sqrt{2\delta}+\sqrt{2}\big)}{\sqrt{\delta}-1}\tag{3.32}
\end{align*}
\section{Some applications II: the $\zeta$ and the $\eta$ functions}
\subsection{The Riemann zeta function}
The Riemann zeta function $\zeta(s)$, proposed by Bernhard Riemann in the year 1859, is an important function associated with the distribution of prime numbers. For $\operatorname{Re}(s)>1$ the Riemann zeta function is defined by the Dirichlet series
\begin{align*}
\zeta(s)=1+\frac{1}{2^s}+\frac{1}{3^s}+\cdots\tag{4.1.1}
\end{align*} 
and is extended to the rest of the complex plane by analytic continuation. In this paper we will not discuss how the Riemann zeta function can be analytically continued to the rest of the complex beyond the $\operatorname{Re}(s)>1$ domain of convergence, but we will be interested in some of the analytically continued results of this function. In the critical strip $0<\operatorname{Re}(s)<1$ and in the region $\operatorname{Re}(s)<0$, the Riemann zeta function satisfies 
\begin{align*}
\zeta(s)=2^s\pi^{s-1}\sin\big(\frac{\pi s}{2}\big)\Gamma(1-s)\zeta(1-s)\tag{4.1.2}
\end{align*}
For $\operatorname{Re}(s)>0$, $\zeta(s)$ and $\eta(s)$ are connected by
\begin{align*}
\eta(s)=\big(1-2^{1-s}\big)\zeta(s)\tag{4.1.3}
\end{align*}
Furthermore, the gamma function $\Gamma(s)$ and $\zeta(s)$ are known to be connected by
\begin{align*}
\Gamma(s)\zeta(s)=\int_{0}^{\infty}{\frac{u^{s-1}}{e^u-1}}du\tag{4.1.4}
\end{align*}
derived from 
\begin{align*}
\Gamma(s)=\int_{0}^{\infty}{e^{-t}t^{s-1}}dt\tag{4.1.5}
\end{align*}
for any $\operatorname{Re}(s)>0$ such that $s\notin\{1, 0, -1, -2, \cdots\}$.  Nice proofs of (4.1.4) can easily be found in the literature such as  in \cite{hv}, but we will restate a proof  here for a purpose. We make a change of variable $t=nu$ to obtain $dt=ndu$. Substituting that into (4.1.5) then gives
\begin{align*}
\Gamma(s)=\int_{0}^{\infty}{e^{-nu}\big(nu\big)^{s-1}}ndu=\int_{0}^{\infty}{n^su^{s-1}e^{-nu}}du\tag{4.1.6}
\end{align*}
Therefore,
\begin{align*}
\Gamma(s)\frac{1}{n^s}&=\int_{0}^{\infty}{u^{s-1}e^{-nu}}du\tag{4.1.7}
\end{align*}
Which follows that
\begin{align*}
\Gamma(s)\sum_{n=1}^{\infty}{\frac{1}{n^s}}=\sum_{n=1}^{\infty}{\int_{0}^{\infty}{u^{s-1}e^{-nu}}du}=\int_{0}^{\infty}{u^{s-1}}\sum_{n=1}^{\infty}{e^{-nu}}du
\end{align*}
Therefore,
\begin{align*}
\Gamma(s)\zeta(s)=\int_{0}^{\infty}{u^{s-1}}\sum_{n=1}^{\infty}{e^{-nu}}du\tag{4.1.8}
\end{align*}
and thus,
\begin{align*}
\Gamma(s)\zeta(s)=\int_{0}^{\infty}{u^{s-1}}\frac{e^{-u}}{1-e^{-u}}du=\int_{0}^{\infty}{\frac{u^{s-1}}{e^u-1}}du\tag{4.1.9}
\end{align*}
We are interested in the geometric series component of (4.1.8). We write
\begin{align*}
\sum_{n=1}^{\infty}{e^{-nu}}=e^{-u}+e^{-2u}+e^{-3u}+\cdots.\tag{4.1.10}
\end{align*}
Applying the geometric series formula in (3.1), we have
\begin{align*}
\sum_{n=1}^{\infty}{e^{-nu}}=e^{-u}\bigg(\frac{e^{-u\delta}-1}{e^{-u}-1}\bigg)=\frac{1-e^{u\delta}}{e^{u\delta}(1-e^u)}
\end{align*}
Therefore,
\begin{align*}
\sum_{n=1}^{\infty}{e^{-nu}}=\frac{e^{-u\delta}-1}{1-e^u}=\frac{1-e^{-u\delta}}{e^u-1}\tag{4.1.11}
\end{align*}
Substituting (4.1.11) back into (4.1.8) gives
\begin{align*}
\Gamma(s)\zeta(s)=\int_{0}^{\infty}{u^{s-1}(\frac{e^{-u\delta}-1}{1-e^u})}du=\int_{0}^{\infty}{\frac{u^{s-1}\big(1-e^{-u\delta}\big)}{e^u-1}}du
\end{align*}
Hence,
\begin{align*}
\Gamma(s)\zeta(s)=\int_{0}^{\infty}{\frac{u^{s-1}\big(1-e^{-u\delta}\big)}{e^u-1}}du\tag{4.1.12}
\end{align*}
One sees that as $n\rightarrow\pm\infty$, the component $1-e^{-u\delta}$ tends to 1 due to the fact that\\ $\displaystyle\lim_{n\to\pm\infty}{e^{-u\delta}}$=$\displaystyle\lim_{n\to\pm\infty}{e^{u\frac{W_n(-\log2)}{\log2}}}\approx0$, and in that case one re-obtains (4.1.9).\par
\subsection{The Euler product formula for $\zeta(s)$}
Riemann's work on what is now named after him, the Riemann zeta function, was built on or motivated by an earlier foundation that Leonhard Euler had laid. Leonhard Euler demonstrated that  for all $\operatorname{Re}(s)>1$, the Dirichlet series in (4.1.1) could be expressed as
\begin{align*}
\zeta(s)=\prod_{p}{(1-p^{-s})^{-1}}=\prod_{p}{(p^s/(p^s-1))}\tag{4.2.1}
\end{align*}
whose unsimplified version is 
\begin{align*}
\zeta(s)=\prod_{p}{\bigg(\sum_{n=0}^{\infty}{p^{-ns}}\bigg)}\tag{4.2.2}
\end{align*}
By comparing (4.2.1) and (4.2.2), we see that
\begin{align*}
\sum_{n=0}^{\infty}{p^{-ns}}=\frac{p^s}{p^s-1}\tag{4.2.3}
\end{align*}
and that is based on the idea that the sum of an infinite geometric series with $|r|<1$ is $\frac{a}{1-r}$. By that approach we know that there is a small error term between $\sum_{n=0}^{\infty}{p^{-ns}}$ and $\frac{p^s}{p^s-1}$. If we were to concern ourselves for a moment, how big or small is this difference? Is this error term significant? In this section we explore answers to these questions. From 
\begin{align*}
\zeta(s)=1+2^{-s}+3^{-s}+4^{-s}+5^{-s}+\cdots
\end{align*}
We arrange the terms such that $(1+2^{-s})$ can be factored out. We obtain

\begin{multline*}
\zeta(s)=\bigg(1+\frac{1}{2^s}\bigg)\bigg(1+\frac{1}{4^s}+\frac{1}{16^s}+\frac{1}{64^s}+\cdots\bigg)\bigg(1+\frac{1}{3^s}+\frac{1}{5^s}+\frac{1}{7^s}+\frac{1}{9^s}+\cdots\bigg)
\end{multline*}
which is expanded and simplified to
\begin{align*}
\zeta(s)=\bigg(\sum_{n=0}^{\infty}{2^{-ns}}\bigg)\bigg(1+\frac{1}{3^s}+\frac{1}{5^s}+\frac{1}{7^s}+\frac{1}{9^s}+\cdots\bigg)\tag{4.2.4}
\end{align*}
We repeat the process by factoring out $(1+3^{-s})$ on the remaining terms to arrive at
\begin{align*}
\zeta(s)=\bigg(\sum_{n=0}^{\infty}{2^{-ns}}\bigg)\bigg(1+\frac{1}{3^s}\bigg)\bigg(1+\frac{1}{3^{2s}}+\frac{1}{3^{4s}}+\frac{1}{3^{6s}}+\cdots\bigg)\bigg(1+\frac{1}{5^s}+\frac{1}{7^s}+\frac{1}{11^s}+\cdots\bigg)
\end{align*}
Which can be simplified to
\begin{align*}
\zeta(s)=\bigg(\sum_{n=0}^{\infty}{2^{-ns}}\bigg)\bigg(\sum_{n=0}^{\infty}{3^{-ns}}\bigg)\bigg(1+\frac{1}{5^s}+\frac{1}{7^s}+\frac{1}{11^s}+\cdots\bigg)\tag{4.2.5}
\end{align*}
Repeating the process infinitely many times over all the primes yields
\begin{align*}
\zeta(s)=\bigg(\sum_{n=0}^{\infty}{2^{-ns}}\bigg)\bigg(\sum_{n=0}^{\infty}{3^{-ns}}\bigg)\bigg(\sum_{n=0}^{\infty}{5^{-ns}}\bigg)\bigg(\sum_{n=0}^{\infty}{7^{-ns}}\bigg)\cdots.\tag{4.2.6}
\end{align*}
Which can further be simplified to
\begin{align*}
\zeta(s)=\prod_{p}\bigg(\sum_{n=0}^{\infty}{p^{-ns}}\bigg)\tag{4.2.7}
\end{align*}
We remember from (2.1.27) that
\begin{align*}
\sum_{n=0}^{\infty}{x^{-n}}=\frac{x}{x-1}-\frac{x^{1-\delta}}{x-1}
\end{align*}
and comparing that to $\sum_{n=0}^{\infty}{p^{-ns}}$ we see that $x=p^s$. This follows, therefore, that
\begin{align*}
\sum_{n=0}^{\infty}{p^{-ns}}=\frac{p^s}{p^s-1}-\frac{p^{s(1-\delta)}}{p^s-1}
\end{align*}
Therefore,
\begin{align*}
\zeta(s)=\prod_{p}{\bigg(\frac{p^s}{p^s-1}-\frac{p^{s(1-\delta)}}{p^s-1}\bigg)}\tag{4.2.8}
\end{align*}
By comparing the Euler's original (4.2.1) with (4.2.8), one sees that the difference $\Delta(p;s)$ is 
\begin{align*}
\Delta(p;s)=-\frac{p^{s(1-\delta)}}{p^s-1}\tag{4.2.9}
\end{align*}
For instance, when $s=2$, the magnitude of the error term in the first prime factor is
\begin{align*}
\Delta(2;2)&=-\frac{2^{2(1-\delta)}}{2^2-1}=-\frac{4^{1-\delta}}{3}=-\frac{4^{1--\frac{W_n(-\log2)}{\log2}}}{3}\\
&=-\frac{4^{1+\frac{W_n(-\log2)}{\log2}}}{3}\tag{4.2.10}
\end{align*}
Some solutions of (4.2.10) are as follows. Note how the solutions decay rapidly to zero as $n\rightarrow\pm\infty$. 
\begin{multline*}
-\frac{4^{1+\frac{W_0(-\log2)}{\log2}}}{3}=0.24074254124109206246977841244522002868847079238710711906...\\ -
0.35029202939416491016186386087571818174009039874188212474... i\\
-\frac{4^{1+\frac{W_{999999999}(-\log2)}{\log2}}}{3}=1.62266893759136173664777462174268821559271968828942... \times 10^{-20}\\ -
1.18423993961989827520959944961560473709484779596843... \times 10^{-28} i
\end{multline*}
The error terms in the other prime factor at $s=2$ are
\begin{align*}
\Delta(3;2)&=-\frac{3^{2(1-\delta)}}{3^2-1}=-\frac{9}{8}\times3^{\frac{W_n(-\log2)}{\log2}}\\
\Delta(5;2)&=-\frac{25}{24}\times25^{\frac{W_n(-\log2)}{\log2}}\\
\Delta(7;2)&=-\frac{49}{48}\times49^{\frac{W_n(-\log2)}{\log2}}\\
\Delta(11;2)&=-\frac{121}{120}\times121^{\frac{W_n(-\log2)}{\log2}}\\
\Delta(13;2)&=-\frac{169}{168}\times169^{\frac{W_n(-\log2)}{\log2}}\\
\Delta(17;2)&=-\frac{289}{288}\times289^{\frac{W_n(-\log2)}{\log2}}
\end{align*}
and many more. As $n\rightarrow\pm\infty$ each of these error terms tends to zero. In conclusion, the  error term in each prime factor is actually very small relative to the actual value of its prime factor, thus $\frac{p^s}{p^s-1}$. In fact, these error terms are basically zero.\par
\subsection{The Dirichlet eta function}
The Dirichlet eta function $\eta(s)$, also called the alternating zeta function, is defined as
\begin{align*}
\eta(s)=\sum_{n=1}^{\infty}{\frac{(-1)^{n-1}}{n^s}}=1-2^{-s}+3^{-s}-4^{-s}+\cdots\tag{4.3.1}
\end{align*}
for all  $\operatorname{Re}(s)>0$. The connection between $\eta(s)$ and $\zeta(s)$ is known to be
\begin{align*}
\eta(s)=\big(1-2^{1-s}\big)\zeta(s)\tag{4.3.2}
\end{align*}
It is known that $\eta(1)=\log2$, but plugging $s=1$ into (4.3.2) yields an undefined result on the right side. So how does that happen? From the idea that
\begin{align*}
\log(1+x)=x-\frac{x^2}{2}+\frac{x^3}{3}-\frac{x^4}{4}+\cdots\tag{4.3.3}
\end{align*}
by substituting $x=1$, one arrives at
\begin{align*}
\log2=1-\frac{1}{2}+\frac{1}{3}-\frac{1}{4}+\frac{1}{5}-\cdots
\end{align*}
which proves that $\eta(1)=\log2$. The question we consider is, how do we explain the observation that $\eta(1)=\log2$ only in the context of (4.3.2)? It turns out (4.3.2) is approximated and in this section we will prove this finding.  Another motivation for dealing with the Dirichlet eta function in this paper is to take advantage of the connection between $\eta(s)$ and $\zeta(s)$ to derive a simpler representation of the harmonic series. In the following theorem, we show if we consider only and only the series in (4.3.1) to derive the relation in (4.3.2), then it is evident that (4.3.2) is approximated. If the small error  between the exact form of (4.3.2) and (4.3.2) itself can be found, then a simpler form of the harmonic series can be derived. We first prove that (4.3.2) is approximated and then give the exact relation for it.
\begin{theorem}
Given that $\eta(s)=1-2^{-s}+3^{-s}-4^{-s}+\cdots$, if $\operatorname{Re}(s)>0$ then
\begin{align*}
\eta(s)=\bigg[\frac{2^{s(\delta+1)}-2(2^{\delta s})+1}{2^s\big[2^{\delta s}-1\big]}\bigg]\zeta(s)\tag{4.3.4}
\end{align*}
\end{theorem}
\begin{proof}
We will first prove that (4.3.2) is approximated. From
\begin{align*}
\eta(s)=1-\frac{1}{2^s}+\frac{1}{3^s}-\frac{1}{4^s}+\frac{1}{5^s}-\cdots\tag{4.3.5}
\end{align*}
Applying the technique of factorization used in subsection 4.2, we have
\begin{align*}
\eta(s)=\bigg(1-\frac{1}{2^s}-\frac{1}{2^{2s}}-\frac{1}{2^{3s}}-\cdots\bigg)\bigg(1+\frac{1}{3^s}+\frac{1}{3^{2s}}+\cdots\bigg)\bigg(1+\frac{1}{5^s}+\frac{1}{5^{2s}}+\cdots\bigg)\cdots\tag{4.3.6}
\end{align*}
\begin{align*}
\eta(s)=\bigg(1-\frac{1}{2^s}\big(1+\frac{1}{2^s}+\frac{1}{2^{2s}}+\cdots\big)\bigg)\bigg(1+\frac{1}{3^s}+\frac{1}{3^{2s}}+\cdots\bigg)\bigg(1+\frac{1}{5^s}+\frac{1}{5^{2s}}+\cdots\bigg)\cdots\tag{4.3.7}
\end{align*}
\begin{align*}
\eta(s)=\bigg(1-\frac{1}{2^s}\big(1+\frac{1}{2^s}+\frac{1}{2^{2s}}+\cdots\big)\bigg)\bigg(\sum_{n=0}^{\infty}{3^{-ns}}\bigg)\bigg(\sum_{n=0}^{\infty}{5^{-ns}}\bigg)\bigg(\sum_{n=0}^{\infty}{7^{-ns}}\bigg)\cdots\tag{4.3.8}
\end{align*}
Which can then be rewritten as
\begin{align*}
\eta(s)=\bigg(\big(1+\frac{1}{2^s}+\frac{1}{2^{2s}}+\cdots\big)^{-1}-\frac{1}{2^s}\bigg)\bigg(\sum_{n=0}^{\infty}{2^{-ns}}\bigg)\bigg(\sum_{n=0}^{\infty}{3^{-ns}}\bigg)\bigg(\sum_{n=0}^{\infty}{5^{-ns}}\bigg)\cdots.\tag{4.3.9}
\end{align*}
But since we have already seen that
\begin{align*}
\zeta(s)=\bigg(\sum_{n=0}^{\infty}{2^{-ns}}\bigg)\bigg(\sum_{n=0}^{\infty}{3^{-ns}}\bigg)\bigg(\sum_{n=0}^{\infty}{5^{-ns}}\bigg)\cdots,
\end{align*}
it follows that (4.3.9) becomes
\begin{align*}
\eta(s)=\bigg(\big(1+\frac{1}{2^s}+\frac{1}{2^{2s}}+\cdots\big)^{-1}-\frac{1}{2^s}\bigg)\zeta(s)\tag{4.3.10}
\end{align*}
Applying the geometric series formula $\frac{a}{1-r}$ on the infinite series in (4.3.10) gives
\begin{align*}
1+\frac{1}{2^s}+\frac{1}{2^{2s}}+\cdots=\frac{1}{1-\frac{1}{2^s}}=\frac{2^s}{2^s-1}\tag{4.3.11}
\end{align*}
Plugging the result of (4.3.11) into (4.3.10)  then gives
\begin{align*}
\eta(s)=\bigg(\bigg(\frac{2^s}{2^s-1}\bigg)^{-1}-\frac{1}{2^s}\bigg)\zeta(s)\tag{4.3.12}
\end{align*}
Which follows that
\begin{align*}
\eta(s)=\bigg(\frac{2^s-1}{2^s}-\frac{1}{2^s}\bigg)\zeta(s)=\bigg(\frac{2^s-2}{2^s}\bigg)\zeta(s)=\bigg(1-2^{1-s}\bigg)\zeta(s)
\end{align*}
Therefore,
\begin{align*}
\eta(s)=\bigg(1-2^{1-s}\bigg)\zeta(s)\tag{4.3.13}
\end{align*}
which is the original relation in (4.3.2). The second approach is as follows. We know that the geometric series in (4.3.11) has $a=1$ and $r=\frac{1}{2^s}$. Plugging these parameters into the geometric series formula in (3.1) gives:
\begin{align*}
1+\frac{1}{2^s}+\frac{1}{2^{2s}}+\cdots=\frac{(\frac{1}{2^s})^{\delta}-1}{\frac{1}{2^s}-1}=\frac{2^{\delta s}-1}{2^{s(\delta-1)}(2^s-1)}\tag{4.3.14}
\end{align*}
We could have equally used the formula $\sum_{n=0}^{\infty}{x^{-n}}=\frac{x}{x-1}-\frac{1}{x^{\delta-1}(x-1)}$ with $x=2^s$ to obtain the same result in (4.3.14). Therefore
substituting the result of (4.3.14) into (4.3.10) gives
\begin{align*}
\eta(s)=\bigg[\bigg(\frac{2^{\delta s}-1}{2^{s(\delta-1)}(2^s-1)}\bigg)^{-1}-\frac{1}{2^s}\bigg]\zeta(s)\tag{4.3.15}
\end{align*}
Which follows that
\begin{align*}
\eta(s)=\bigg[\frac{2^{s(\delta-1)}(2^s-1)}{2^{\delta s}-1}-\frac{1}{2^s}\bigg]\zeta(s)
\end{align*}
Which is further written as
\begin{align*}
\eta(s)=\bigg[\frac{2^{\delta s}\big(2^s-1\big)-2^{\delta s}+1}{2^s\big[2^{\delta s}-1\big]}\bigg]\zeta(s)=\bigg[\frac{2^{s(\delta+1)}-2^{\delta s}-2^{\delta s}+1}{2^s\big[2^{\delta s}-1\big]}\bigg]\zeta(s)
\end{align*}
Therefore,
\begin{align*}
\eta(s)=\bigg[\frac{2^{s(\delta+1)}-2(2^{\delta s})+1}{2^s\big[2^{\delta s}-1\big]}\bigg]\zeta(s)\tag{4.3.16}
\end{align*}
as required.
\end{proof}
One sees that with a little manipulation of  (4.3.16) one arrives at
\begin{align*}
\eta(s)=\bigg[\frac{2^{s(\delta+1)}-2(2^{\delta s})+1}{2^{s(\delta+1)}-2^s}\bigg]\zeta(s)=\bigg[\frac{1-2^{1-s}+2^{-s(\delta+1)}}{1-2^{-\delta s}}\bigg]\zeta(s)
\end{align*}
From which one further sees that the terms $2^{-s(\delta+1)}$ and $2^{-\delta s}$ decay very rapidly to zero as $n\rightarrow\pm\infty$ since $\delta$ depends on $W_n$ and $\operatorname{Re}(s)>0$. In that case, one then arrives at the already known relation which states that
\begin{align*}
\eta(s)=\big(1-2^{1-s}\big)\zeta(s)
\end{align*}
as already seen in (4.3.2). But for some important reasons which we will later see, we will still go with the idea that
\begin{align*}
\eta(s)=\bigg[\frac{2^{s(\delta+1)}-2(2^{\delta s})+1}{2^s\big[2^{\delta s}-1\big]}\bigg]\zeta(s)
\end{align*}
By that same approach we could show that for  $\operatorname{Re}(s)<0$,
\begin{align*}
\eta(-s)=\bigg[\frac{2^{s+1}-2^{s(\delta+1)}-1}{2^{\delta s}-1}\bigg]\zeta(-s)\tag{4.3.17}
\end{align*}
but we will not concern ourselves about  proving it in this paper.
\subsubsection{\Large The harmonic series}
The harmonic series, usually regarded as the value of $\zeta(s)$ at $s=1$, is an important infinite series in mathematics due to its connection with the prime numbers. Euler's proof of the infinitude of prime numbers suggests that the behavior of the harmonic series is somehow related with the distribution of prime numbers, \cite{notes}. Furthermore, Euler's proof of primes' infinitude involves the idea that it is the summation of the reciprocals of prime numbers that causes the harmonic series to diverge and that is possible because there are infinitely many prime numbers.\par
Speaking of a proof of the divergence of the harmonic series, there are numerous approaches in the literature for doing so. Some of these approaches involve the idea of comparing the behavior of the harmonic series with a known behavior of a related infinite series and then making deductions. Such approaches can be hard to follow sometimes. In this section we will give two new representations  of the harmonic series and prove using those methods that the harmonic series diverges. The harmonic series is commonly represented as
\begin{align*}
\zeta(1)=\sum_{k=1}^{\infty}{\frac{1}{k}}=1+\frac{1}{2}+\frac{1}{3}+\frac{1}{4}+\cdots=\infty.\tag{4.3.1.1}
\end{align*}
From (4.3.16), we remember that
\begin{align*}
\eta(s)=\bigg[\frac{2^{s(\delta+1)}-2(2^{\delta s})+1}{2^s\big[2^{\delta s}-1\big]}\bigg]\zeta(s)
\end{align*}
Plugging $s=1$ into the above relation gives
\begin{align*}
\eta(1)=\bigg[\frac{2^{(\delta+1)}-2(2^{\delta })+1}{2\big[2^{\delta }-1\big]}\bigg]\zeta(1)
\end{align*}
\begin{align*}
\log2&=\bigg[\frac{2^{\delta+1}-2^{\delta+1 }+1}{2\big[2^{\delta }-1\big]}\bigg]\zeta(1)\\
\log2&=\bigg[\frac{1}{2\big[2^{\delta }-1\big]}\bigg]\zeta(1)
\end{align*}
Therefore,
\begin{align*}
\zeta(1)=2\big[2^{\delta}-1\big]\log2\tag{4.3.1.2}
\end{align*}
But we remember that $\delta=2^\delta=-\frac{W_n(-\log2)}{\log2}$. On that note, we have
\begin{align*}
\zeta(1)=2\bigg[-\frac{W_n(-\log2)}{\log2}-1\bigg]\log2=-2\big[\log2+W_n(-\log2)\big]
\end{align*}
Therefore the harmonic series, denoted by $\zeta(1)$, is given by
\begin{align*}
\zeta(1)=-2\big[\log2+W_n(-\log2)\big]\tag{4.3.1.3}
\end{align*}
We show some of the solutions of (4.3.1.3) below. Note how the real and the imaginary parts of the solutions grow bigger and bigger toward $\tilde\infty$ as $n\rightarrow\pm\infty$. This, as we have seen from previous cases in section 3, shows that the harmonic series indeed diverges. 
\begin{multline*}
-2\big[\log2+W_0(-\log2)\big]\\=-0.2430471428665835983605365717644769258177599837351707030...\\ -
2.172922314730940848930566050905504325083853368632426901... i\\
-2\big[\log2+W_{-1}(-\log2)\big]\\=-0.2430471428665835983605365717644769258177599837351707030...\\ +
2.172922314730940848930566050905504325083853368632426901... i\\
-2\big[\log2+W_1(-\log2)\big]\\=3.48685848447467640681800979173824201426011212844133009055...\\ -
15.0829564115551999194410614719764865099606395615389014483... i\\
-2\big[\log2+W_{-2}(-\log2)\big]\\=3.48685848447467640681800979173824201426011212844133009055...\\+
15.0829564115551999194410614719764865099606395615389014483... i\\
-2\big[\log2+W_{10000}(-\log2)\big]\\=21.4432163834098854395227571741079321210490350858570638860...\\ -
125666.847372911457288756710995905214075491077681109327037... i
\end{multline*}
\begin{multline*}
-2\big[\log2+W_{-10001}(-\log2)\big]\\=21.4432163834098854395227571741079321210490350858570638860...\\ +
125666.847372911457288756710995905214075491077681109327037... i\\
-2\big[\log2+W_{999999999}(-\log2)\big]\\=44.46901728525495132622711474304605111793813593691208...\\ -
1.256637060493439498578309429671659338741364104305815... \times 10^{10} i\\
-2\big[\log2+W_{-1000000000}(-\log2)\big]\\=44.46901728525495132622711474304605111793813593691208...\\ +
1.256637060493439498578309429671659338741364104305815... \times10^{10} i
\end{multline*}
By looking at the sizes of the real parts of these solutions of the harmonic series relative to the sizes of $n$ as $n\rightarrow\pm\infty$, one sees that the harmonic series grows very slowly.\par
The second representation of the harmonic series uses the relation in (2.1.21). That,
\begin{align*}
1+\frac{m}{2}+\frac{m^2}{3}+\cdots=\frac{m^{\delta-1}}{\delta}+\frac{m^{\delta-2}}{\delta-1}+\frac{m^{\delta-3}}{\delta-2}+\cdots
\end{align*}
When $m=1$, we obtain
\begin{align*}
1+\frac{1}{2}+\frac{1}{3}+\frac{1}{4}+\cdots=\frac{1}{\delta}+\frac{1}{\delta-1}+\frac{1}{\delta-2}+\frac{1}{\delta-3}+\cdots\tag{4.3.1.4}
\end{align*}
But $\delta=-\frac{W_n(-\log2)}{\log2}$. Therefore
\begin{align*}
1+\frac{1}{2}+\frac{1}{3}+\frac{1}{4}+\cdots&=-\frac{\log2}{W_n(-\log2)}-\frac{\log2}{W_n(-\log2)+\log2}-\frac{\log2}{W_n(-\log2)+2\log2}\\
&-\frac{\log2}{W_n(-\log2)+3\log2}-\frac{\log2}{W_n(-\log2)+4\log2}-\cdots\tag{4.3.1.5}
\end{align*}
But we remember from (4.3.1.3) that $\zeta(1)=-2\big[\log2+W_n(-\log2)\big]$, and substituting that into (4.3.1.5) gives
\begin{align*}
-2\big[\log2+W_n(-\log2)\big]&=-\frac{\log2}{W_n(-\log2)}-\frac{\log2}{W_n(-\log2)+\log2}-\frac{\log2}{W_n(-\log2)+2\log2}\\
&-\frac{\log2}{W_n(-\log2)+3\log2}-\frac{\log2}{W_n(-\log2)+4\log2}-\cdots
\end{align*}
\begin{align*}
\frac{-2\big[\log2+W_n(-\log2)\big]}{-\log2}&=\frac{1}{W_n(-\log2)}+\frac{1}{W_n(-\log2)+\log2}+\frac{1}{W_n(-\log2)+2\log2}\\
&+\frac{1}{W_n(-\log2)+3\log2}+\frac{1}{W_n(-\log2)+4\log2}+\cdots
\end{align*}
\begin{align*}
2\bigg[1+\frac{W_n(-\log2)}{\log2}\bigg]&=\frac{1}{W_n(-\log2)}+\frac{1}{W_n(-\log2)+\log2}+\frac{1}{W_n(-\log2)+2\log2}\\
&+\frac{1}{W_n(-\log2)+3\log2}+\frac{1}{W_n(-\log2)+4\log2}+\cdots
\end{align*}
But we also remember from (2.10) that 
\begin{align*}
\phi=1+2+2^2+2^3+\cdots=-1-\frac{W_n(-\log2)}{\log2}
\end{align*}
and on that note we have
\begin{align*}
-2\big[1+2+2^2+2^3+\cdots\big]&=2+\frac{2W_n(-\log2)}{\log2}\\&=\frac{1}{W_n(-\log2)}+\frac{1}{W_n(-\log2)+\log2}+\frac{1}{W_n(-\log2)+2\log2}\\
&+\frac{1}{W_n(-\log2)+3\log2}+\frac{1}{W_n(-\log2)+4\log2}+\cdots\tag{4.3.1.6}
\end{align*}
Some solutions of (4.3.1.6) are computed below. Note how the real and the imaginary parts of these solutions grow bigger and bigger toward $\tilde\infty$ as $n\rightarrow\pm\infty$. This is a sign of divergence. One should not confuse (4.3.1.6) with the harmonic series representation in (4.3.1.5). The fact that (4.3.1.6) diverges means the harmonic series diverges, as we see their connection from (4.3.1.5).
\begin{multline*}
2+\frac{2W_0(-\log2)}{\log2}=0.35064290771585155537187081123679352005078515963673834662...\\ +
3.1348642476992957221171487823859738550666151484084351263... i\\
2+\frac{2W_{-1}(-\log2)}{\log2}=0.35064290771585155537187081123679352005078515963673834662... \\-
3.1348642476992957221171487823859738550666151484084351263... i\\
2+\frac{2W_{9999}(-\log2)}{\log2}=-30.935733390839434708251937393886716439092103299954771114...\\ +
181280.80806843448062030568837327049614764711199663830658... i\\
2+\frac{2W_{-10000}(-\log2)}{\log2}=-30.935733390839434708251937393886716439092103299954771114...\\ -
181280.80806843448062030568837327049614764711199663830658... i
\end{multline*}
The result below is due to Euler, as he was the first to deeply study it. Using the idea that 
\begin{align*}
\zeta(s)=\bigg(\sum_{n=0}^{\infty}{2^{-ns}}\bigg)\bigg(\sum_{n=0}^{\infty}{3^{-ns}}\bigg)\bigg(\sum_{n=0}^{\infty}{5^{-ns}}\bigg)\cdots
\end{align*}
 and substituting $s=1$ gives
\begin{align*}
\zeta(1)=\bigg({\sum_{n=0}^{\infty}{2^{-n}}}\bigg)\bigg({\sum_{n=0}^{\infty}{3^{-n}}}\bigg)\bigg({\sum_{n=0}^{\infty}{5^{-n}}}\bigg)\bigg({\sum_{n=0}^{\infty}{7^{-n}}}\bigg)\cdots
\end{align*}
Applying the geometric series formula $\Pi=a\big(\frac{r^{\delta}-1}{r-1}\big)$ on each factor of the  infinite product gives
\begin{align*}
\zeta(1)=\bigg(\frac{2^{-\delta}-1}{2^{-1}-1}\bigg)\bigg(\frac{3^{-\delta}-1}{3^{-1}-1}\bigg)\bigg(\frac{5^{-\delta}-1}{5^{-1}-1}\bigg)\bigg(\frac{7^{-\delta}-1}{7^{-1}-1}\bigg)\cdots
\end{align*}
which follows by simplification that
\begin{align*}
\zeta(1)=2\big[1-2^{-\delta}\big]\times\frac{3}{2}\big[1-3^{-\delta}\big]\times\frac{5}{4}\big[1-5^{-\delta}\big]\times\frac{7}{6}\big[1-7^{-\delta}\big]\times\cdots.
\end{align*}
By rearrangement of terms we obtain
\begin{align*}
\zeta(1)=\bigg(2\times\frac{3}{2}\times\frac{5}{4}\times\frac{7}{6}\times\frac{11}{10}\times\cdots\bigg)\bigg(1-2^{-\delta}\bigg)\bigg(1-3^{-\delta}\bigg)\bigg(1-5^{-\delta}\bigg)\bigg(1-7^{-\delta}\bigg)\cdots.
\end{align*}
From which we then see that
\begin{align*}
\zeta(1)\zeta(\delta)=\frac{2\times3\times5\times7\times11\times13\times\cdots}{1\times2\times4\times6\times10\times12\times\cdots}
\end{align*}
We have already seen that $\zeta(1)=-2\big[\log2+W_n(-\log2)\big]$ and $\delta=-W_n(-\log2)/\log2$. Therefore,
\begin{align*}
-2\big[\log2+W_n(-\log2)\big]\zeta(-W_n(-\log2)/\log2)=\frac{2\times3\times5\times7\times11\times13\times\cdots}{1\times2\times4\times6\times10\times12\times\cdots}
\end{align*}
But we see that as $n\rightarrow\pm\infty$ the term $\zeta(-W_n(-\log2)/\log2)\rightarrow1$ and that makes sense since we know that $\delta$ is actually $\delta=1+1+2+2^2+2^3+\cdots$. In that case one arrives at
\begin{align*}
-2\big[\log2+W_n(-\log2)\big]=\frac{2\times3\times5\times7\times11\times13\times\cdots}{1\times2\times4\times6\times10\times12\times\cdots}\tag{4.3.1.7}
\end{align*}
which was discovered by Leonhard Euler in the eighteenth century. The difference is that Euler did not know that $\zeta(1)=-2\big[\log2+W_n(-\log2)\big]$.\par
In [2], García and Marco prove that the super-regularized product over all the prime numbers is $4\pi^2$. If we assume accuracy of their result, then from (4.3.1.7) we obtain
\begin{align*}
1\times2\times4\times6\times10\times12\times16\times\cdots&=-\frac{4\pi^2}{2\big[\log2+W_n(-\log2)\big]}\\
\prod_{p_n}{(p_n-1)}&=-\frac{2\pi^2}{\log2+W_n(-\log2)}\tag{4.3.1.8}
\end{align*}
From the result of (4.3.1.8) we see that as $n\rightarrow\pm\infty$, $\displaystyle-\frac{2\pi^2}{\log2+W_n(-\log2)}\rightarrow0$, which, in fact, must be the case in order  for (4.3.1.7) to diverge. Also, multiplying both sides of (4.3.1.8) by $\log2$ gives
\begin{align*}
\log2\prod_{p_n}{(p_n-1)}=-\frac{2\pi^2\log2}{\log2+W_n(-\log2)}=2\pi^2\bigg(\sum_{n=0}^{\infty}{2^n}\bigg)^{-1}
\end{align*}
since we know from (2.10) that $\phi=1+2+2^2+\cdots=-1-W_n(-\log2)/\log2$. Therefore,
\begin{align*}
2\pi^2=\bigg(\log2\sum_{n=0}^{\infty}{2^n}\bigg)\bigg(\prod_{p_n}{(p_n-1)}\bigg)\tag{4.3.1.9}
\end{align*}
\section{Some unrigorous results}
Let us sidestep traditions and pretend for a moment that there was no such a function as the Riemann zeta function. Let us also assume that since the Riemann zeta function does not exist, we have no idea of its so-called trivial zeros which occur at all the  negative even $n$. Then let us define some function, say $\varrho(s)$, by
\begin{align*}
\varrho(s)=1+2^s+3^s+4^s+5^s+\cdots\tag{5.1}
\end{align*}
By the Euler product approach, we could rewrite (5.1) as
\begin{align*}
\varrho(s)=\bigg(\sum_{n=0}^{\infty}{2^{ns}}\bigg)\bigg(\sum_{n=0}^{\infty}{3^{ns}}\bigg)\bigg(\sum_{n=0}^{\infty}{5^{ns}}\bigg)\bigg(\sum_{n=0}^{\infty}{7^{ns}}\bigg)\cdots\tag{5.2}
\end{align*}
But we remember from (2.1.26) that $\sum_{n=0}^{\infty}{x^n}=x^{\delta-1}\sum_{n=0}^{\infty}{x^{-n}}$. By this approach, we could transform (5.2) into the following result:
\begin{align*}
\varrho(s)&=\bigg(2^{s(\delta-1)}\sum_{n=0}^{\infty}{2^{-ns}}\bigg)\bigg(3^{s(\delta-1)}\sum_{n=0}^{\infty}{3^{-ns}}\bigg)\bigg(5^{s(\delta-1)}\sum_{n=0}^{\infty}{5^{-ns}}\bigg)\cdots\\
&=\big(2\times3\times5\times7\times\cdots\big)^{s(\delta-1)}\bigg(\sum_{n=0}^{\infty}{2^{-ns}}\bigg)\bigg(\sum_{n=0}^{\infty}{3^{-ns}}\bigg)\bigg(\sum_{n=0}^{\infty}{5^{-ns}}\bigg)\cdots\tag{5.3}
\end{align*}
Based on (5.2), we see that
\begin{align*}
\varrho(-s)=\bigg(\sum_{n=0}^{\infty}{2^{-ns}}\bigg)\bigg(\sum_{n=0}^{\infty}{3^{-ns}}\bigg)\bigg(\sum_{n=0}^{\infty}{5^{-ns}}\bigg)\cdots
\end{align*}
So (5.3) becomes
\begin{align*}
\varrho(s)=\varrho(-s)\big(2\times3\times5\times7\times11\times\cdots\big)^{s(\delta-1)}\tag{5.4}
\end{align*}
From (5.4) one sees that in terms of the almighty Riemann zeta function, we have
\begin{align*}
\zeta(-s)=\zeta(s)\big(2\times3\times5\times7\times11\times\cdots\big)^{s(\delta-1)}\tag{5.5}
\end{align*}
which will not hold at any even $n$ since we have $\zeta(-s)$ and $\zeta(s)$ on opposite sides of the equation and that is due to the presence of the trivial zeros of $\zeta(s)$. From that formula we further see that the infinite product over the primes will vary with values of $s$ and that  could certainly be telling us something interesting about these prime numbers. Therefore, since the infinite product over all $p$ varies with $s$ in this case, it is not possible to replace the infinite product over the primes  in (5.5) with García and Marco's result which states that the super-regularized product over all the prime numbers is $4\pi^2$. Let us just assume for a moment that the relation in (5.5) is actually correct. Then we observe the following:
when $s=489$, we have
\begin{align*}
\zeta(-489)=\big[2\times3\times5\times7\times\cdots\big]^{489(\delta-1)}\zeta(489)
\end{align*}
which gives,
\begin{align*}
\big[2\times3\times5\times7\times\cdots\big]^{489(\delta-1)}=-4.352545655837123484889282539871247099567386436677185257\\465395744156381927210480072700330299261823650910387... \times 10^{713}
\end{align*}
Therefore,
\begin{align*}
\big[2\times3\times5\times7\times\cdots\big]^{(\delta-1)}=28.798838375431513979727792227935021915438295336551280656...\\ +
0.18502139847178457723424976648813747104332649080311856456... i\\
\end{align*}
When $s=491$, we have
\begin{align*}
\big[2\times3\times5\times7\times\cdots\big]^{(\delta-1)}=\sqrt[491]{\frac{\zeta(-491)}{\zeta(491)}}=28.916650510146632196520893211219095990702485367634\cdots
\end{align*}
And when $s=507$, we have
\begin{align*}
\big[2\times3\times5\times7\times\cdots\big]^{(\delta-1)}=\sqrt[507]{\frac{\zeta(-507)}{\zeta(507)}}=29.85437601678774667123912113393667379536402578\cdots
\end{align*}
When $s=601$, we have
\begin{align*}
\big[2\times3\times5\times7\times\cdots\big]^{(\delta-1)}=\sqrt[601]{\frac{\zeta(-601)}{\zeta(601)}}=35.36250344681743069025439843179082329124629134417...\\ +
0.1848512361343942233986287432135067584710057... i
\end{align*}
When $s=2019$, we have 
\begin{align*}
\big[2\times3\times5\times7\times\cdots\big]^{(\delta-1)}=\sqrt[2019]{\frac{\zeta(-2019)}{\zeta(2019)}}=118.42187068137957389115884655189084428470907521...
\end{align*}
When $s=3000+2i$, we have
\begin{align*}
\big[2\times3\times5\times7\times\cdots\big]^{(\delta-1)}=\sqrt[3000+2i]{\frac{\zeta(-3000-2i)}{\zeta(3000+2i)}}=176.012331418729098119946889468358867117...\\ -
0.71232444285334541154087526334826173578... i
\end{align*}
When $s=5000+2i$, we have
\begin{align*}
\big[2\times3\times5\times7\times\cdots\big]^{(\delta-1)}=\sqrt[5000+2i]{\frac{\zeta(-5000-2i)}{\zeta(5000+2i)}}=293.1275991242533831525164044517...-\\
0.711689211684005594786759086472672... i
\end{align*}
The value of the term $\big[2\times3\times5\times7\times\cdots\big]^{(\delta-1)}$ actually varies with $s$. What could that be telling us about the prime numbers associated with it? Interestingly, some values such as $s=491, 507, 2019, \cdots$ give only real values to this term. Why is that so?\par
\section{Concluding Remarks} In this paper we present an approach which helps us to deal with convergent and divergent cases of all infinite geometric series with $r\neq1$. By this same approach we are able to gain insight into how to approach some other cases of non-geometric infinite series and that includes the harmonic series, leading to two new proofs that the harmonic series diverges.\par
We have seen that all the divergent infinite series considered in this paper have infinitely many complex solutions, which indicate that such series indeed do not converge. Apart from using these numerous complex solutions to show that they are signs of divergence, in what other ways could one apply these solutions to make an impact?\\

$\textbf{Acknowledgement}$. My sincere gratitudes go to Dr. Kofi Adanu, a researcher at the Alabama Transportation Institute, Dr. Larry Gratton and Dr. Jay Baltisberger, Professors at  Berea College, for their advice and motivation throughout this work.

\small Cletus Bijalam Mbalida\\
101 Chestnut Street\\
College Post Office Box 968\\
Berea, Kentucky 40403-1516\\
The United States of America\\
$\textit{Email:}$ \href{mailto:mbalidac@berea.edu}{mbalidac@berea.edu}
\end{document}